\documentclass[11pt,oneside]{article}   	
\usepackage{amsmath,amsthm,graphicx,amssymb,enumerate,color,natbib,bbm,rotating,tikz,authblk,romannum,array,enumitem,dsfont,booktabs,multirow,mathtools,float}
\usepackage{caption}
\usepackage{dutchcal}
\usepackage[pagewise]{lineno}
\usepackage[font=small]{caption}
\usepackage[mathscr]{euscript}
\usepackage{makecell}
\usepackage[top=1in, bottom=1in, left=1in, right=1in]{geometry}
\usepackage{subcaption}
\usepackage{float}
\usepackage{booktabs}
\usepackage{siunitx}
\usetikzlibrary{chains,shapes.multipart,shapes,calc,fit}
\usepackage[percent]{overpic}
\usetikzlibrary{automata,positioning}
\usepackage[bottom]{footmisc}
\usepackage[hidelinks=true]{hyperref}
\usepackage{setspace}
\usepackage[T1]{fontenc}
\usepackage[utf8]{inputenc}
\usepackage[para,online,flushleft]{threeparttable}
\usepackage{algorithm}
\usepackage{algpseudocode}
\usepackage[all]{xy}
\usepackage{array}  
\setcitestyle{round}
\DeclareMathAlphabet{\mathpzc}{OT1}{pzc}{m}{it}

\newtheoremstyle{ModifiedStyle}
{\topsep} 
{3pt} 
{} 
{} 
{\bfseries} 
{.} 
{.5em} 
{} 
\theoremstyle{ModifiedStyle}

\newtheorem{proposition}{Proposition}

\newtheorem{assumption}{Assumption}

\newtheorem{to-do}{To-Do}

\newcounter{subroutine}
\newcounter{mainalgorithm}
\newenvironment{Subroutine}[1][htb]
  {\refstepcounter{subroutine}%
   \begin{algorithm}[#1]}
  {\end{algorithm}}

\newenvironment{MainAlgorithm}[1][htb]
  {\refstepcounter{mainalgorithm}%
   \begin{algorithm}[#1]}
  {\end{algorithm}}

\newcommand{\BSmall}[1]{\mkern-1.7mu\raisebox{-1.2pt}{\scalebox{0.9}{$\scriptscriptstyle B$}}}
\newcommand{\NSmall}[1]{\mkern-1.7mu\raisebox{-1.2pt}{\scalebox{0.9}{$\scriptscriptstyle N$}}}

\DeclareMathOperator*{\maximize}{maximize}

\newcommand{\nquad}{\kern-1em}

\makeatletter
\newcommand{\Rom}[1]{\expandafter\@slowromancap\romannumeral #1@}
\newcommand{\Biggg}{\bBigg@{2.5}}
\newcommand{\vast}{\bBigg@{3}}
\newcommand{\Vast}{\bBigg@{3.5}}
\newcommand{\massive}{\bBigg@{4.5}}
\newcommand{\Massive}{\bBigg@{6}}
\newcommand{\thickhline}{%
	\noalign {\ifnum 0=`}\fi \hrule height 1pt
	\futurelet \reserved@a \@xhline}
\newcolumntype{"}{@{\hskip\tabcolsep\vrule width 1pt\hskip\tabcolsep}}

\newcommand{\PreserveBackslash}[1]{\let\temp=\\#1\let\\=\temp}
\newcolumntype{C}[1]{>{\PreserveBackslash\centering}p{#1}}
\newcolumntype{R}[1]{>{\PreserveBackslash\raggedleft}p{#1}}
\newcolumntype{L}[1]{>{\PreserveBackslash\raggedright}p{#1}}
\makeatother

\title{{\textbf{Dynamic control of stochastic matching systems in
heavy traffic: An effective computational method for high-dimensional problems}}}

\author[$\star$]{Barı\c{s} Ata}
\author[$\star$]{Yaosheng Xu}
\affil[$\star$]{{\small Booth School of Business, The University of Chicago}}
\date{{\small \vspace{-7mm} \today}}

\makeatletter
\let\LN@align\align
\let\LN@endalign\endalign
\renewcommand{\align}{\linenomath\LN@align}
\renewcommand{\endalign}{\LN@endalign\endlinenomath}
\let\LN@gather\gather
\let\LN@endgather\endgather
\renewcommand{\gather}{\linenomath\LN@gather}
\renewcommand{\endgather}{\LN@endgather\endlinenomath}
\makeatother

\begin{document}
	
	\vspace*{-22mm}
	
	{\let\newpage\relax\maketitle}
	
	
	
	\vspace*{-9mm}
	\begin{abstract}
Bipartite matching systems arise in many settings where agents or tasks from two distinct sets must be paired dynamically under compatibility constraints. We consider a high-dimensional bipartite matching system under uncertainty and seek an effective dynamic control policy that maximizes the expected discounted total value generated by the matches minus the congestion-related costs. To derive a tractable approximation, we focus attention on balanced, high-volume systems, i.e., the heavy-traffic regime, and derive an approximating Brownian control problem. We then develop a computational method that relies on deep neural network technology for solving this problem. To show the effectiveness of the policy derived from our computational method, we compare it to the benchmark policies available in the extant literature in the context of the original matching problem. In the test problems attempted thus far, our proposed policy outperforms the benchmarks, and its derivation is computationally feasible for dimensions up to 100 or more.

	\end{abstract}
	
	
	\pagenumbering{arabic}
	\doublespacing
	
	\setlength{\abovedisplayskip}{8pt}
	\setlength{\belowdisplayskip}{8pt}
	\setlength{\abovedisplayshortskip}{8pt}
	\setlength{\belowdisplayshortskip}{8pt}
	
\section{Introduction}
Dynamic matching problems arise in many operations management applications, such as ride-hailing platforms, gig marketplaces, paired kidney exchange, dating markets, public housing allocation, etc. In such systems, jobs of different types--e.g., agents and tasks--arrive at the system over time and must be matched subject to compatibility constraints. In this paper, we focus on a bipartite matching system in which a system manager makes matching decisions dynamically over time under uncertainty.

In our formulation, jobs from different classes arrive to the system according to a Poisson process. When a job arrives, if there is a compatible match already waiting in the system, the system manager may match them and earn the associated value from the match. Upon matching, both jobs leave the system immediately. Otherwise, the arriving job joins a buffer to wait for a suitable match. As jobs wait in the system, the system manager incurs a holding cost per waiting job. Additionally, jobs may abandon as they wait in the system, resulting in abandonment penalties.

The system manager seeks a dynamic matching policy to maximize the expected discounted value generated by the matches, performed over time minus the holding and abandonment costs. We model this system as a queueing network, or a stochastic processing network, as described in \cite{Harr2003}. Crucially, we seek an effective solution to this problem for the high-dimensional case. Solving such problems is intractable beyond a few dimensions using standard dynamic programming methods due to the curse of dimensionality, which is the case even for the smallest test problem we consider; see Section \ref{sec:testproblems} for the X model.

To derive a tractable formulation, we first approximate the queueing network formulation by a Brownian control problem (BCP), focusing on the balanced, large-flow asymptotic regime, i.e., the heavy-traffic regime. More specifically, our approximating problem involves the drift-rate control of a reflected Brownian motion (RBM). In particular, we derive a drift-rate control problem as the approximation, as opposed to a singular control problem as is typically done in the heavy-traffic literature (see, for example, \cite{HarrWein1989}). Next, to solve the resulting drift control problem, we adopt the deep learning-based computational method proposed by \cite{AtaHarrSi2024}. We test the effectiveness of the resulting control policy on a variety of test problems by comparing it with benchmark policies drawn from the literature. It is worth mentioning that all comparisons are made in the context of the original queueing formulation of the matching problem, \emph{not} in the context of the approximating BCP, even though our proposed policy is derived from the BCP.

To position the contribution of our paper in the literature, we briefly mention some of the related prior research; see Section~2 for a detailed review of the relevant literature. Several researchers have used first-order approximations, i.e., approximations based on deterministic fluid models, to derive effective matching policies; see, for example, \cite{GurvWard2014}. We also start with a nominal plan of long-run matching rates derived from a linear program that uses only job arrival rates and the value generated from each possible match. However, we further refine the policy prescription dynamically based on the system state for which the solution of the approximating BCP is used crucially. We observe that doing so significantly improves the test performance in our test problems. We also contribute to the analysis of stochastic processing networks by illustrating the applicability of the modeling framework proposed by \cite{Harr2003} and \cite{AtaHarrSi2025} to a dynamic matching problem, and providing an effective control policy for the high-dimensional setting by adopting the computational framework developed in \cite{AtaHarrSi2024}.

We test the performance of the proposed policy on nine test problems. We begin with a simple $2 \times 2$ matching model, where each class can be matched with one of two classes. We refer to this as the "X" model, reflecting the network’s shape. We consider three instances of the X model that are differentiated by their abandonment rates. Next, we consider more complex networks, starting with the $4 \times 4$ example (of dimension 8) from \cite{KerAshGur2024}. We also consider two variants of that model. We refer to these three $4\times 4$ models as the "Zigzag" network models due to their graph structure and use them as building blocks for constructing larger test problems. In particular, we consider two $12 \times 12$ test problems (of dimension 24). The first problem concatenates the three variants of the Zigzag network model by introducing additional basic activities to connect them. In contrast, the second problem concatenates three copies of the most complex Zigzag model using nonbasic activities; see Section 4 for definitions of basic versus nonbasic activities. Lastly, we consider a $60 \times 60$ (120-dimensional) test problem that is constructed by concatenating copies of the second 24-dimensional model. For all test problems considered, even the X model, characterizing the optimal policy via standard dynamic programming techniques is not computationally tractable. Thus, we rely on benchmark policies derived in prior literature. The proposed policy significantly outperforms the benchmark policies that we consider.

The rest of the paper is structured as follows. Section \ref{sec:literature} reviews related literature. Section \ref{sec:model} introduces the model. In Section \ref{sec:aproxbcp}, we derive the approximating Brownian control problem in the heavy-traffic asymptotic regime. Then we characterize the optimal drift-rate control policy given the optimal value function obtained using our computational method. Section \ref{sec:compute} outlines the computational approach used to solve the Brownian control problem. Section \ref{sec:policy} introduces our proposed policy. Benchmark policies drawn from the literature are presented in Section \ref{sec:benchmark}. We present the test problems in Section \ref{sec:testproblems} and the computational results in Section \ref{sec:results}. Section \ref{sec:conclusion} concludes. Appendices \ref{app:deriveBCP} through \ref{app:NNarchi} provide the derivation of the Brownian control problem, proofs, model parameters for the high-dimensional test problems, and the details of the neural network architecture. 

\section{Literature review}\label{sec:literature}
\cite{GalSha1962} started a large literature on matching. Their work was motivated by college admissions. The authors introduced the solution concept of an optimal stable match and illustrated it with their deferred acceptance procedure. Since then, matching models have been used to study a variety of applications such as school choice (e.g., \cite{AbduSonm2003}, \cite{AshlShi2015}, and \cite{FeigKanoLoSeth2020}); kidney exchange (e.g., \cite{RothSonmUnve2004}, \cite{RothSonmUnverDelm2006}); public housing allocation (e.g., \cite{Kaplan1986} and \cite{ArnoShi2020}); ride-hailing (e.g., \cite{WanZhaZha2023}); financial engineering (e.g., \cite{MagMoaZhe2021}) and dating markets (e.g., \cite{KanoSaba2021}). See \cite{Koji2017}, \cite{ChadEecSmi2017} and \cite{HuanTanWajc2024} for surveys of this literature. 

Much of the literature on matching, especially in the economics literature, models incentives of different parties, whereas we focus solely on the resource allocation problem, or the first-best decisions, see \cite{LaffTiro1993}. Additionally, most of the extant literature focuses on static models, while we focus on a dynamic matching problem. As such, the most closely related papers to ours can naturally be modeled as (two-sided) queueing systems. We divide those papers into two groups. Roughly speaking, the first group focuses on performance evaluation whereas the second group focuses on optimization. 

\cite{AdanWeis2014} considers a bipartite matching system under the first-come-first-served, assign-longest-idle-server (FCFS-ALIS) policy and derives a product-form expression for the stationary distribution of this system when the service capacity is sufficient. \cite{AdanWeis2014} extends the earlier work by \cite{CalKapWei2009}, \cite{AdanWeis2012}  and \cite{VissAdanWeis2012}.
\cite{DiaBar2019} considers a two-sided queueing system, where each side has patient and impatient customers. Additionally, customers belong to either the high- or low-priority class. Within each priority class, customers are matched according to the FCFS policy. The authors derive closed-form expressions for the steady-state distributions of the high- and low-priority queues. They also derive closed-form expressions for other steady-state performance measures. 

\cite{BusGupMai2013} studies the stability of bipartite matching systems under certain policies. The authors show that the stability region for the match-the-longest (ML) policy is maximal for any bipartite graph. Relatedly, \cite{AfehCaldGupt2022} studies the design of a bipartite matching topology for multiclass, multiserver queueing systems operating under the FCFS-ALIS discipline. The objective is to balance two competing goals in a Pareto sense: minimizing congestion costs and maximizing rewards from matches; also see \cite{HillaCaldGupt2024}. \cite{BukeChen2017} and \cite{LiuGonKul2015} consider fluid and diffusion limits of matching models where either a match occurs with a certain probability (\cite{BukeChen2017}) or everyone matches when there is an available match (\cite{LiuGonKul2015}). In both cases, the system state can be reduced to a one-dimensional representation. 

Among the papers that consider control of matching systems, perhaps the most closely related one is \cite{HuZhou2022}. The authors consider a finite-horizon, discrete-time matching problem that is similar to ours at a high level. They formulate the problem as a dynamic program and characterize its structural properties. The authors characterize sufficient conditions under which a certain hierarchical priority policy is optimal that combines a first stage greedy policy with a second-stage threshold policy at each decision epoch. \cite{LeeLiuLiuZhang2021} considers a two-sided queueing system as a model of a time-varying production system. Adopting a fluid model, the authors optimize the (nonstationary) production rate over a finite horizon. \cite{BlanReimShahWeinWu2022} considers a parsimonious two-sided matching system and studies a threshold policy which uses a single threshold to decide when the matches should be made. The authors characterize the asymptotically optimal threshold policies under different distributional assumptions; also see \cite{MertNaxPrad2024}.

\cite{GurvWard2014} seeks to minimize cumulative holding cost over a finite horizon. The authors adopt a discrete-review framework; see \cite{Magla1998} and \cite{AtaKuma2005}. They show that a myopic discrete-review matching policy is optimal on the fluid scale. \cite{NazaStol2019} uses a greedy primal-dual approach to derive an asymptotically optimal policy that seeks to maximize the revenue rate while maintaining stable queues. \cite{OzkWar2020} uses a matching model to study a time-varying ride-hailing system. The authors use a fluid-approximation and derive an asymptotically optimal policy in the fluid scale, utilizing the associated continuous linear program; also see \cite{WanZhaZha2023} for a similar study. 

\cite{VarBumMagWan2023} considers a problem of dynamic pricing and matching for a bipartite matching system. Using a fluid-model based pricing policy coupled with the max-weight matching policy, the authors establish an asymptotically optimal performance. They derive further refinements of their policy with improved suboptimality guarantees in an asymptotic sense. \cite{AveDeLeviWard2024} considers a bipartite matching model. The authors derive a discrete-review matching policy using the associated fluid model and establish its asymptotic optimality on the fluid scale.  There is also a stream of work in online bipartite matching in an adversarial setting \cite{KarpVaziVazi1990} where agents do not wait if they are not matched immediately; also see \cite{KerAshGur2024, KerAshGur2025}, \cite{WeiXuYu2023} and the references therein for similar models. Some studies explore how random graph structures impact matching outcomes, such as phase transitions or recoverability thresholds, see \cite{AndeItaiGamaKano2017}, \cite{MohaMoorXu2021}, \cite{MaoWuXuYu2025}, etc. These streams of literature are peripheral to our work. 

The discrete-flow bipartite matching model developed in Section 3 is intractable. To derive a tractable approximation, we consider it in the balanced, high-volume, i.e., the heavy traffic regime. \cite{harr1998, harr2000, Harr2003} pioneered an approach to approximate such dynamic control problems with the so-called Brownian control problems, which in turn can be reduced to singular control problems. Additionally, recent work by \cite{AtaHarrSi2025} develops a novel modification of \cite{Harr2003} that approximates the original scheduling problem by a drift-rate control of reflected Brownian motion (RBM). We follow the latter approach to derive our approximating Brownian control problem in Section \ref{sec:aproxbcp}.

Drift-rate control problems have been used in the literature to study control problems arising in various applications. \cite{AtaHarrShep2005} considers a one-dimensional drift-rate control problem motivated by a power control problem in wireless communication. \cite{MatogVand2011} considers a drift rate control problem where a system manager incurs a fixed cost to change the drift rate. \cite{GhoshWeer2007, GhoshWeer2010} extend \cite{AtaHarrShep2005}. \cite{Ata2006}, \cite{CeliMagl2008} and \cite{AtaBar2023} build on \cite{AtaHarrShep2005} to address different problems arising in make-to-order production systems. \cite{RubiAta2009} uses a drift-rate control formulation to study a scheduling problem for a parallel-server queueing system with abandonments. \cite{BarIMarPerr2007} uses a drift-rate control formulation for a macroeconomics application. \cite{AtaLeeSonm2019, AtaTongLeeFiel2024} study volunteer scheduling problems using drift-rate control formulations; also see \cite{AtaLeeTong2024}. \cite{AlwaAtaZhou2024} uses a drift-rate control formulation to study a problem motivated by a ride-hailing application.

 The aforementioned papers all study one-dimensional control problems. Indeed, drift-rate control problems in higher dimensions do not admit closed-form solutions. Neither are their solutions tractable computationally beyond a few dimensions using historical grid-based methods, e.g., finite-element methods. However, recent papers \cite{AtaHarrSi2024}, \cite{AtaKas2025} and \cite{AtaZhou2025} solve drift-rate control problems arising in different contexts in dimensions tens to hundreds utilizing the deep BSDE methods that will be reviewed next. 

In order to solve the drift-rate control problems, one considers the associated HJB equation that is a semilinear PDE. Historically, the state-of-the-art method for solving them has been the finite-element method, which is a grid-based method. Thus, it suffers from the curse of dimensionality. In the last decade, there has been a breakthrough in solving partial differential equations. Broadly speaking, two primary approaches have been developed: Physics-Informed Neural Networks (PINNs) and methods based on Backward Stochastic Differential Equations (BSDEs).
PINNs, proposed by \cite{RaisPerdKarn2019}, provide a meshless method for solving partial differential equations (PDEs) by integrating observational data with fundamental physical laws. In this approach, the target function is approximated by a neural network, whose gradients and Hessians are obtained via automatic differentiation. The neural network parameters are optimized by minimizing a loss function that incorporates the PDE residuals along with boundary or initial condition errors. A key difficulty with this method is the high computational cost associated with evaluating Hessians using automatic differentiation, particularly as PDE dimensionality grows. Nonetheless, several strategies exist to accelerate these computations, such as those discussed by \cite{HeLiShiGao2023} and \cite{HuWangZhou2024}.

The deep BSDE approach, initially developed by \cite{EHanJent2017, HanJenE2018}, exploits the connection between semilinear PDEs and backward stochastic differential equations (BSDEs), as originally established by \cite{PardPeng1990} and \cite{MaProtYong1994}. Specifically, \cite{HanJenE2018} approximates both the (initial) value function and its spatial gradient using neural networks across discrete time steps, ultimately approximating the terminal value function. The neural network parameters are optimized by minimizing a loss function defined as the discrepancy between this approximation and the given terminal condition. This method offers two notable advantages relevant to our work: it circumvents the computationally intensive Hessian evaluation and facilitates efficient sampling from the state space for training purposes when a reference policy is available. See \cite{EHanJen2022}, \cite{BeckHutzJentKuck2023}, \cite{ChesKawaiShinoToshi2023}, \cite{HanJentE2025} for reviews of this literature. 

\section{Model}\label{sec:model}
\begin{figure}[H]
		\centering	\includegraphics[width=0.35\textwidth]{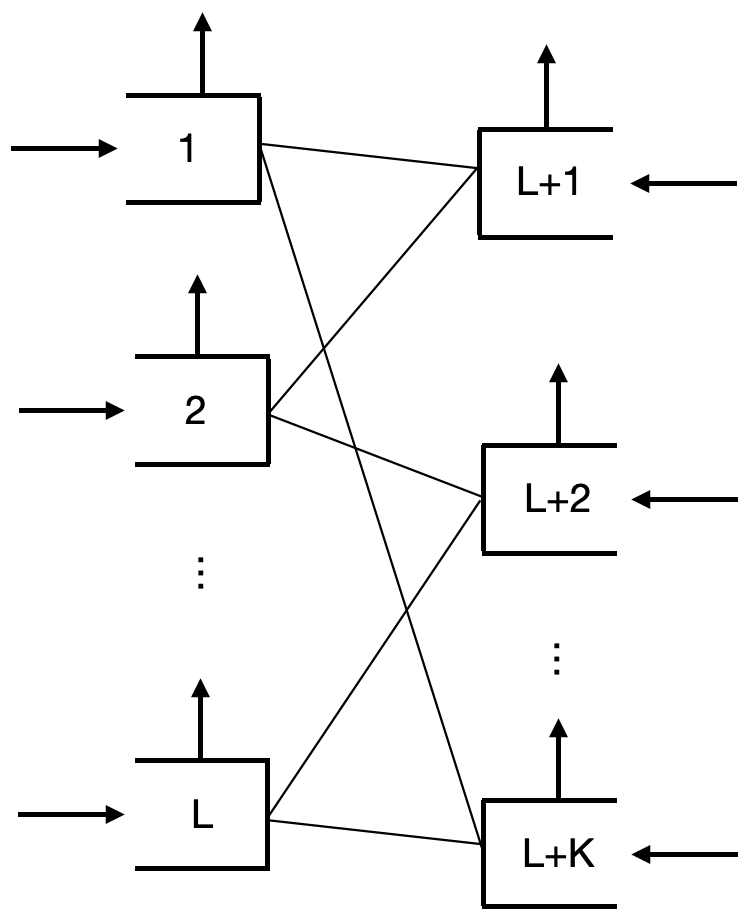}
		\caption{A bipartite matching model.}
        \label{fig:match}
\end{figure}

We consider the bipartite matching model displayed in Figure \ref{fig:match} that comprises two distinct sets of classes. The first set (on the left in Figure \ref{fig:match}) consists of $L$ classes, whereas the second one (on the right in Figure \ref{fig:match}) has $K$ classes, i.e., there are $I = L + K$ classes in total. Each class has an associated buffer in which jobs of that class waiting to be matched reside. In what follows, we use the terms buffer and class interchangeably. The first set of buffers is denoted by $\mathcal{L}$ (mnemonic for left), whereas the second one is denoted by $\mathcal{R}$ (mnemonic for right). We index the buffers so that $\mathcal{L} = \{1, \ldots, L\}$ and $\mathcal{R} = \{L+1, \ldots, L+K\}$.
There are $J$ matching activities. Activity $j$ corresponds to matching a job in buffer $\ell(j) \in \mathcal{L}$ with a job in buffer $r(j) \in \mathcal{R}$ for $j = 1, \ldots, J$. The activities are denoted by the lines connecting the buffers in Figure \ref{fig:match}. To avoid trivial complications,  we assume for every buffer $i = 1, \ldots, I$ that there exists an activity that can match it with another buffer. Associated with each activity is a value $v_j > 0$ for $j = 1, \ldots, J$.

We let $E_i(t)$ denote the cumulative number of class $i$ jobs arrived by time $t$, and assume the processes $\{E_i(t), t \geq 0\}$ for $i = 1, \ldots, I$ are independent Poisson processes with rate $\lambda_i > 0$. We also let $E(t) = (E_i(t))$ denote the corresponding vector-valued arrival process and $\lambda = (\lambda_i)$ is the $I$-dimensional vector of arrival rates. The Poisson assumption simplifies the exposition, but it is not essential for our analysis. All that is needed is that the process $\{E(t), t \geq 0\}$ satisfies a functional central limit theorem.

When a job arrives, if a matching activity is feasible, i.e., both buffers involved in the match are nonempty, then the system manager may perform the match and earn the associated value from it. Otherwise, the job waits in the designated buffer. Jobs within each class are homogeneous and are matched in the order they arrive. That is, we assume first-come-first-served (FCFS) for matching within each class. As they wait, jobs in class $i$ incur a holding cost rate of $h_i > 0$ for $i = 1, \ldots, I$. Additionally, as they wait in the queue, jobs may abandon. We model the time-to-abandon as an exponential random variable with rate $\gamma_i$ for class $i$ jobs. For each abandoning class $i$ job, the system manager incurs a cost of $a_i$ for $i = 1, \ldots, I$.

Letting $Q_i(t)$ denote the number of class $i$ jobs in the system at time $t$, the $I$-dimensional queue-length process $Q(t) = (Q_i(t))$ represents the system state. The cumulative number of class $i$ jobs abandoned up to time $t$, denoted by $A_i(t)$, is modeled as follows:
\begin{align}\label{eq:aban}
A_i(t)=N_i\left(\int_0^t \gamma_i Q_i(s) d s\right), \quad t\geq 0,\quad i=1, \ldots, I, \end{align}
where $N_i(\cdot)$ is a Poisson process with rate one, and we assume $N_1, \ldots, N_I$ and $E_1, \ldots, E_I$ are mutually independent.

The system manager decides which matching activities to undertake over time dynamically. We let $T_j(t)$ denote the cumulative number of matches made via activity $j$, i.e., the total number of matches that involve buffers $\ell(j)$ and $r(j)$, until time $t$. Then the $J$-dimensional process $T(t) = (T_j(t))$ for $t \geq 0$ denotes the system manager's control. We assume that matches happen instantaneously and require that 
\begin{equation}\label{eq:T}
T(\cdot) \text{ is non-anticipating and non-decreasing with } T(0) = 0.
\end{equation}
Additional restrictions for admissibility of a control $T(t)$ will be specified below.

In order to facilitate our description of the evolution of the system state mathematically, we introduce the $I \times J$ dimensional matching matrix $R$, where its $j$th column, denoted by $R^j$, corresponds to activity $j$. We let
\begin{equation}\label{eq:R}
R^j = e_{\ell(j)} + e_{r(j)}, \quad j = 1, \ldots, J,
\end{equation}
where $e_i$ denotes the $I$-dimensional unit vector that has 1 in the $i$th position and has zeros elsewhere. In words, $R^j$ indicates the number of jobs deleted from each buffer per unit of activity $j$. Naturally, activity $j$ removes a job from each of the buffers $\ell(j)$ and $r(j)$; no other buffer is affected. The $j$th column of the matching matrix $R$ captures this transition of the system state. Then, given a control $T(\cdot)$, the evolution of the system state can be described as follows:
\begin{equation}\label{eq:Q_dyn}
Q(t) = Q(0) + E(t) - R T(t) - A(t), \quad t \geq 0.
\end{equation}
For a control $T(\cdot)$ to be feasible, we require the resulting queue length vector to be nonnegative. That is,
\begin{equation} \label{eq:Q_pos}
Q(t) \geq 0, \quad t \geq 0.
\end{equation}
Under a feasible matching control $T(\cdot)$, the system manager earns the cumulative net value of $\Pi(t)$ until time $t$, where
\begin{equation}\label{eq:Pi}
\Pi(t) = v \cdot T(t) - \int_0^t h \cdot Q(s) \, ds - a \cdot A(t), \quad t \geq 0. 
\end{equation}
The system manager seeks a feasible matching control $T(\cdot)$ so as to
\begin{align}
\text{maximize} \quad \mathbb{E} \left[ \int_0^\infty e^{-r t} \, d\Pi(t) \right]\quad \text{ subject to } (\ref{eq:aban}) - (\ref{eq:Pi}),
\end{align}
where $r$ denotes the interest rate for discounting. 

This problem is intractable. Thus, building on \cite{Harr2003} and \cite{AtaHarrSi2025}, we derive an approximating Brownian control problem in the heavy traffic limiting regime, i.e., the large and balanced volume regime.  We solve it by building on the approach developed in \cite{AtaHarrSi2024}, which in turn builds on the seminal work of \cite{HanJenE2018}; also see \cite{AtaKas2025}.

\section{Approximating Brownian Control Problem}\label{sec:aproxbcp}

Following \cite{Harr2003}, we first introduce a static planning problem that helps the system manager choose the nominal long-run activity rates. 

\noindent\textbf{Static planning problem.} Choose a $J$-dimensional column vector $x$ so as to 
	\begin{equation}\label{eq:spp}
	\begin{aligned}
	&\maximize\quad v\cdot x \quad \text{subject to }\quad Rx=\lambda, \quad x\geq 0. 
	\end{aligned}
	\end{equation}
The static planning problem considers an idealized setting, where one seeks to maximize the total reward rate $v \cdot x$ associated with the activity rate vector $x$. The constraint $Rx = \lambda$ ensures the activity rate vector $x$ matches all arriving jobs, whereas the constraint $x \geq 0$ captures the natural requirement that the activity rates are nonnegative.

Our heavy traffic assumption entails Assumptions 1 and 2 below.

	\begin{assumption}\label{as:bf}(\textit{Balanced flow assumption}). There exists a unique optimal solution $x^*$ to the static planning problem (\ref{eq:spp}) such that
$Rx^*=\lambda$.
\end{assumption}
It is worth emphasizing that unlike the settings of \cite{Harr2003} or \cite{AtaHarrSi2024}, our model involves no capacity constraints because we assume the matches happen instantaneously. Nonetheless,
we build on their framework to derive the approximating Brownian control problem, making changes as needed.

As in Harrison (2003), we call activity $j$ \textit{basic} if $x_j^* > 0$, and call it \textit{nonbasic} if $x_j^* = 0$. We let $b$ denote the number of basic activities and assume without loss of generality that activities $1, \ldots, b$ are the basic ones. Also, we write \begin{align}\label{eq:os}
	x^* = \begin{bmatrix}
	x_B^*\\ x_N^* 
	\end{bmatrix}, 
	\end{align}
	where $x_B^*=(x_1^*,\ldots,x_b^*)'$ is the $b$-dimensional vector of basic activity rates, whereas $x_N^*$ is the $J-b$ dimensional vector of zeros. Correspondingly, we partition the matching matrix $R$ too so that 
    \begin{align}R= [H\quad  N],\end{align}
	where $H$ is an $I\times b$ dimensional matrix; its columns correspond to the basic activities whereas $N$ is an $I\times (J-b)$ dimensional matrix and its columns correspond to the nonbasic activities. To facilitate future analysis, we also partition the value rates as
$$
v = \begin{bmatrix} v_B \\ v_N \end{bmatrix},
$$
where $v_B = (v_1, \ldots, v_b)' $ and $v_N = (v_{b+1}, \ldots, v_J)' $.

    As is usually done in the heavy traffic literature, we consider a sequence of systems indexed by $n = 1, 2, \ldots$ to derive the approximating Brownian control problem. The superscript $n$ will be attached to the quantities of interest associated with the $n$th system. We focus on the asymptotic regime where the system grows with $n$ as specified in Assumption \ref{as:lf}.

    \begin{assumption}\label{as:lf}(\textit{Large flow assumption}). We assume that $\lambda^n =n\lambda$ for $n\geq 1$ where $\lambda$ is an $I$-dimensional vector of strictly positive entries. 
    \end{assumption}
It is worth mentioning that for every buffer there must exist a basic activity matching it with another buffer, because we assume $\lambda_i > 0$ for all $i=1,\ldots,I$ and that $Rx^* = \lambda$. As mentioned earlier, our heavy traffic assumption comprises of Assumptions \ref{as:bf} and \ref{as:lf}. It is evident from them that if one ignores randomness, thus the holding and abandonment costs, the system manager would choose the instantaneous activity rate vector $n x^*$ for the $n$th system. However, under  uncertainty, she may benefit from small adjustments to this activity rate vector dynamically depending on the system state. In addition, one expects the queue lengths to be of second order relative to the system size, i.e., order $\sqrt{n}$. Thus, we define their scaled version as follows:
\begin{equation}\label{eq:zn}
Z^n(t) = \frac{Q^n(t)}{\sqrt{n}}, \quad t \geq 0. 
\end{equation}
To facilitate the derivation of the Brownian control problem, we define a drift-rate function $\theta: \mathbb{R}_+^I \to [-\eta, \eta]^J$, where $\eta > 0$ will be viewed as a tuning parameter below. We also introduce the $J$-dimensional non-decreasing process $Y(t)$. As in \cite{AtaHarrSi2024b}, the system manager’s control is the drift rate function $\theta(\cdot)$ in the approximating Brownian control problem, whereas the process $Y(\cdot)$ will serve as a backup control to ensure the (scaled) queue length process remains in the state space. More will be said about the process $Y(\cdot)$ below.

We partition both $\theta(\cdot)$ and $Y(\cdot)$ with respect to the indices of basic versus nonbasic activities as follows:
\begin{align}
	\theta (z) = \begin{bmatrix}
	\theta_B(z)\\ \theta_N(z) 
	\end{bmatrix}, z\geq0,\text{ and  }  Y(t) = \begin{bmatrix}
	Y_B(t)\\ Y_N(t)
	\end{bmatrix}, t\geq 0.
\end{align}
Here, $\theta_B(t)$ contains the first $b$ components of $\theta(t)$ and $\theta_N(t)$ contains the last $J-b$ components. Similarly, $Y_B(t)$ corresponds to the first $b$ components of $Y(t)$. Crucially, we set $Y_N(t) = 0$ for all $t \geq 0$.

As argued in \cite{Harr2003}, any policy worthy of consideration should satisfy $T^n(t) \approx n x^* t$ for $t \geq 0$ when $n$ is large, but the system manager can choose second-order, i.e., order $\sqrt{n}$, deviations from this nominal plan. To capture this, our approximation assumes that the system manager’s control $T^n(\cdot)$ in the $n$th system satisfies the following for $t\geq 0$, 
	\begin{align}
	T_j^n (t)& = nx_j^* t -\sqrt{n} \int_{0}^t\theta_j\left(
    Z^n(s)\right)ds -\sqrt{n} Y_j(t) + o(\sqrt{n}),\quad j=1,\ldots, b,\label{eq:tb}\\
	T_j^n(t) &= -\sqrt{n} \int_{0}^t \theta_j\left(
    Z^n(s)\right)ds + o(\sqrt{n}), \quad j=b+1,\ldots, J.\label{eq:tnb}
	\end{align}
 We also require that 
    \begin{align}
    &Y_j(t) \text{ is non-anticipating, non-decreasing with } Y_{j}(0)=0,\quad j=1,\ldots,b,\label{eq:y}\\
	&Y_j(t) \text{ can only increases at time } t \text{ when either } Z^n_{\ell(j)}(t)=0 \text{ or } Z^n_{r(j)}(t)=0, \quad j=1,\ldots,b.  \label{eq:yb}\end{align}
The requirement in Equation (\ref{eq:yb}) can be expressed equivalently as follows:
\begin{align}\label{eq:complem}\int_{0}^{\infty}Z^n_{\ell(j)}(t)Z^n_{r(j)}(t) dY_j(t)=0,\quad  j=1,\ldots,b. \end{align}
In order to ensure the control $T^n$ is non-decreasing, we further require that 
	\begin{align}
	\theta_B(z)\leq \sqrt{n} x_B^*\quad \text{and} \quad\theta_N(z)\leq 0 \quad\text{for } z\geq 0.
	\end{align}
Also, we require by definition that 
\begin{align}
|\theta_j(z)|\leq \eta,\quad j=1,\ldots, J.
\end{align}
As mentioned earlier, we view $\eta$ as a tuning parameter.
Lastly, we define the centered and scaled cumulative cost process
\begin{align}\label{eq:xi_n}
\xi^n(t) = \frac{(v \cdot x^*) n t - \Pi^n(t)}{\sqrt{n}}, \quad t \geq 0,
\end{align}
where for $i=1,\ldots, I$ and $t \geq 0$,
$$\Pi^n(t) = v \cdot T^n(t) - \int_0^t h \cdot Q^n(s) \, ds - a \cdot A^n(t),\text{ and }A_i^n(t) = N_i(\int_0^t \gamma_i Q_i^n(s) d s).$$ In a similar fashion to Harrison (2003), Appendix \ref{app:deriveBCP} derives the approximating Brownian control problem formally as $n$ gets large. In the approximating Brownian control problem, we fix a large system parameter $n$ and replace the performance-relevant processes $\xi^n, Z^n$ with their formal limits $\xi, Z$, which satisfy the following: for $t\geq 0$, 
\begin{align}
	\xi(t)
	&= \int_{0}^{t}c \cdot Z(s) ds + \int_{0}^{t} v\cdot \theta(Z(s))ds + \int_{0}^{t}v_B\cdot d Y_B(s),\label{eq:xi}\\
	Z(t) &= z + X(t) +\int_{0}^{t}R\,\theta(Z(s))ds + HY_B(t)- \int_0^t \Gamma Z(s) d s,\text{ where } 	\label{eq:z}  \\
    \theta(z)& \in \Theta := \left\{ \theta : \theta_B \leq \sqrt{n} x_B^*, \ \theta_N \leq 0,\ |\theta_j| \leq \eta, \ j = 1, \ldots, J \right\},\label{eq:theta}
	\end{align}
where $X$ is an $I$-dimensional driftless Brownian motion with a diagonal covariance matrix $\Lambda = \mathrm{diag}(\lambda_1, \ldots, \lambda_I)$, and $\Gamma = \mathrm{diag}(\gamma_1,\ldots,\gamma_I)$ is the diagonal matrix of abandonment rates for various classes. Also, $c_i = h_i + \gamma_i a_i$ is the effective holding cost rate for buffer $i$, which combines the holding cost rate $h_i$ and the expected abandonment cost rate $\gamma_i a_i$ for each class $i$ job; and $c = (c_i)$ is the $I$-dimensional vector of effective holding cost rates. Lastly, $\eta$ is the tuning parameter, and $z$ is the initial system state that we incorporate to facilitate future analysis.

We call a drift-rate control function $\theta(\cdot): \mathbb{R}_+^I \to [-\eta, \eta]^J$ admissible if the constraints $(\ref{eq:y})-(\ref{eq:yb})$ and $(\ref{eq:xi}) - (\ref{eq:theta})$ hold. Given an admissible drift-rate function $\theta(\cdot)$, its discounted cost starting in state $z$, denoted by $J(z; \theta)$, is given as follows:
\begin{equation}\label{eq:vf_b}
J(z; \theta) = \mathbb{E}_z^{\theta} \left[ \int_0^{\infty} e^{-rt} \, d\xi(t) \right],
\end{equation}
where $\mathbb{E}_z$ denotes the conditional expectation starting in state $z$. One can equivalently write 
\begin{equation}\label{eq:vf}
J(z; \theta) = \mathbb{E}_z^{\theta}  \left[ \int_0^\infty e^{-r t} \left( c \cdot Z(t) + v \cdot \theta\left(Z(t)\right) \right) dt + \int_0^\infty e^{-r t} v_B \cdot dY_B(t) \right].
\end{equation}
We now define the optimal value function for the Brownian control problem as follows:
$$
V(z) = \inf J(z; \theta), \quad z \in \mathbb{R}_+^I,
$$
where the infimum is taken over admissible drift rate controls. 

The next subsection introduces the Hamilton-Jacobi-Bellman (HJB) equation and the optimal policy for the Brownian control problem that is derived from the HJB equation.

\subsection{HJB equation and the optimal policy for the Brownian control problem}

To solve the Brownian control problem, we next consider the associated HJB equation. As a preliminary, let $f: \mathbb{R}_+^I\rightarrow \mathbb{R}$ be an arbitrary $C^2$ (twice continuously differentiable) function, and let $\nabla f$ denote its gradient vector. We also define the second-order differential operator $\mathcal{L}$ via
$$
\mathcal{L}f = \frac{1}{2} \sum_{i=1}^I \lambda_i \frac{\partial^2 f}{\partial z_i^2},
$$
and a first-order differential operator $\mathcal{D} = (\mathcal{D}_1, \ldots, \mathcal{D}_b)'$ via $\mathcal{D}f = H'\nabla f$. Then, by a similar argument as in \cite{AtaHarrSi2024}, the HJB equation can be stated as follows:
\begin{align}\label{eq:hjb}
\mathcal{L} V(z) + \min_{\theta \in \Theta} \left\{ (R'\nabla V(z) + v ) \cdot \theta \right\}+ (c -  \Gamma \nabla V(z)) \cdot z = r V(z), \quad z \in \mathbb{R}_+^I,
\end{align}
with boundary conditions
\begin{align}\label{eq:hjbb}
\mathcal{D}_j V(z) = -v_j \quad \text{if either } z_{\ell(j)}=0 \text{ or } z_{r(j)} = 0, \quad j = 1, \ldots, b.
\end{align}

\paragraph{Optimal drift-rate control policy for the Brownian control problem.}
The optimal policy for the Brownian control problem is given by
$$
\theta^*(z) = \arg\min_{\theta \in \Theta} \left\{ (R'\nabla V(z) + v) \cdot \theta\right\} , \quad z \in \mathbb{R}_+^I.
$$
To characterize this policy further, we let $U(z) := R' \nabla V(z)+ v$ for $z \in \mathbb{R}_+^I$, and note that for $j = 1, \ldots, J$,
\begin{align}\label{eq:u}
U_j(z) = \frac{\partial V(z)}{\partial z_{\ell(j)}} + \frac{\partial V(z)}{\partial z_{r(j)}} + v_j, \quad z \in \mathbb{R}_+^I.
\end{align}
Then, the optimal policy $\theta^*(\cdot)$ is given explicitly as follows: for $j=1,\ldots, J$, 
\begin{equation}\label{eq:theta*}
\theta_j^*(z) =
\begin{cases}
  -\eta & \text{if } U_j(z) \geq 0, \\
  \min\{\sqrt{n} x_j^*, \eta\} & \text{if } U_j(z) < 0.
\end{cases}
\end{equation}
Thus, once $\nabla V(z)$ is computed, we can obtain the optimal policy (\ref{eq:theta*}) as the solution to the Brownian control problem. The computation of $\nabla V(z)$ will be discussed in Section \ref{sec:compute}.

The index function $U_j(z)$ includes the immediate value $v_j$ generated by activity $j$ as well as the change in the value function through its partial derivatives that are associated with activity $j$. The partial derivative terms capture the future effect of the current state change. In that sense, it is reminiscent of the shadow price concept in linear programming. Thus, we refer to $U_j(z)$ as the \emph{shadow value} of activity $j$ in state $z$.

\section{Computational Method} \label{sec:compute}
Our computational method builds on that of \cite{AtaHarrSi2024}. We begin by specifying a \textit{reference policy} to generate the sample paths of the system state to be visited during neural network training. Roughly speaking, we seek to choose a reference policy so that the resulting paths tend to occupy the parts of the state space that we expect the optimal policy to visit often. To be specific, our reference policy chooses a constant drift rate $\theta(z) = \tilde{\theta} \in \Theta$ in every state $z \in \mathbb{R}_+^I$. The corresponding reference process, denoted by $\tilde{Z}(t)$, satisfies the following:
$$
\tilde{Z}(t) = z + X(t) + R \tilde{\theta} t + H \tilde{Y}_B(t) - \int_0^t \Gamma \tilde{Z}(s) \, ds, \quad t \geq 0.
$$

\begin{Subroutine}[h]
    \captionsetup{labelformat=empty}
		\caption{\textbf{Subroutine 1:} Euler discretization scheme}
		\label{alg:euler}
        \textbf{Input:} Drift vector $\tilde{\theta}$, covariance matrix $\Lambda$, value vector $v_B$, reflection matrix $R$, time horizon $T$, step size $h$ (assume $N = T/h$ is integer), and initial state $\tilde{Z}(0) = z$.
		
        \textbf{Output:} Discretized RBM path $\tilde{Z}$, boundary-pushing increments $\Delta Y$, and Brownian increments $\delta$
		\begin{algorithmic}[1]
			\Function{Discretize}{$T$, $h$, $z$}
			\State Partition $[0, T]$ into $N$ intervals of size $h$
			\State Generate $N$ i.i.d. Gaussian  vectors with mean $0$ and covariance $h\Lambda$: $\delta_0, \ldots, \delta_{N-1}$
			\For{$k \gets 0$ to $N-1$}
			\State $x \gets \tilde{Z}(kh) + \delta_k + R\tilde{\theta}h -\Gamma \tilde{Z}(kh) h$
			\State $\tilde{Z}((k+1)h), \Delta Y(kh) \gets \textsc{Skorokhod}(x;v_B)$ \Comment{See Subroutine~\ref{alg:skorokhod}}
			\EndFor
			\State \textbf{return} $\tilde{Z}(h), \tilde{Z}(2h), \ldots, \tilde{Z}(Nh); \Delta Y(0), \Delta Y(h), \ldots, \Delta Y((N-1)h);$ and $\delta_0, \ldots, \delta_{N-1}$
			\EndFunction
		\end{algorithmic}
\end{Subroutine}
Next, we provide a key identity that guides our definition of the loss function to be minimized during neural network training. For the key identity (\ref{eq:sde1}) below, let
$$
F(z, x) = - (R \tilde{\theta}) \cdot x + c \cdot z + \min_{\theta \in \Theta} \left\{ (x' R + v)  \cdot \theta\right\}\quad \text{for } z\geq 0,x\in\mathbb{R}^I.
$$
\begin{proposition}\label{equivSDE}
		If $V: \mathbb{R}_{+}^I \rightarrow \mathbb{R}$ is a $C^2$ function with polynomial growth that satisfies the $HJB$ equation  (\ref{eq:hjb}) and (\ref{eq:hjbb}), then it also satisfies the following identity almost surely for any $T>0$:
		\begin{equation}\label{eq:sde1}
		\begin{aligned}
		e^{-r T} V(\tilde{Z}(T))-V(\tilde{Z}(0)) & =\int_0^T e^{-r t} \nabla V(\tilde{Z}(t)) \cdot \mathrm{d} W(t) \\
		& -\int_0^T e^{-r t} v
        _B \cdot \mathrm{d} \tilde{Y}_{B}(t)-\int_0^T e^{-r t} F(\tilde{Z}(t), \nabla V(\tilde{Z}(t))) \mathrm{d} t. 
		\end{aligned}
		\end{equation}
	\end{proposition}
\noindent Proposition \ref{equivSDE} provides the motivation for the loss function we use in neural network training; see Appendix \ref{app:equivSDE} for its proof. Additionally, one can prove a converse to Proposition \ref{equivSDE} as done in \cite{AtaHarrSi2024} with minor modifications, which we omit for brevity.
Following \cite{AtaHarrSi2024}, we approximate the value function $V(\cdot)$ and its gradient $\nabla V(\cdot)$ by the deep neural networks $\hat{V}(\cdot; w_1)$ and $G(\cdot; w_2)$, respectively, with associated parameter vectors $w_1$ and $w_2$. To find an approximate solution of the stochastic equation (\ref{eq:sde1}), we define the loss function
\begin{equation}\label{eq:loss}\begin{aligned}
	\ell\left(w_1, w_2\right)= & \mathbb{E}\left[\left(e^{-r T} V_{w_1}(\tilde{Z}(T))-V_{w_1}(\tilde{Z}(0))+\int_0^T e^{-r t} v_B \cdot \mathrm{d} \tilde{Y}_{B}(t)\right.\right. \\
	& \left.\left.-\int_0^T e^{-r t} G_{w_2}(\tilde{Z}(t)) \cdot \mathrm{d} W(t)+\int_0^T e^{-r t} F\left(\tilde{Z}(t), G_{w_2}(\tilde{Z}(t))\right) \mathrm{d} t\right)^2\right].
	\end{aligned}\end{equation}

\noindent In Equation (\ref{eq:loss}), one can choose the initial state $\tilde{Z}(t)$ randomly and the expectation is calculated with respect to the sample path distribution of the reference process $\tilde{Z}$; see Subroutine \ref{alg:compute} for details. As mentioned in Step 4 of the algorithm, the process is run continuously. That is, the terminal state of one iteration becomes the initial state of the next one, as an approximation to starting the reference process with its steady-state distribution.

\begin{Subroutine}[h]
\captionsetup{labelformat=empty} 
    \caption{\textbf{Subroutine 2:} Solve the Skorokhod Problem}\label{alg:skorokhod}
    
    \textbf{Input:} State $x \in \mathbb{R}^J$, reflection matrix $H$, value vector $v_B $, and basic sets $J(i) = \{j=1,\ldots, b:\ell(j) = i \text{ or } r(j) = i\}$, for $i=1,\ldots, I$. 
    
    \textbf{Output:} A solution to the Skorokhod problem $z, y \in \mathbb{R}_+^J$
    
    \begin{algorithmic}[1]
        \Function{Skorokhod}{$x; v_B$}
        \State Set  tolerance $\epsilon = 10^{-8}$
        \State Initialize $z \gets x$
        \While{there exists $i$ such that $z_i < -\epsilon$}
        \State $B \gets\{i : z_i < \epsilon\}$ \Comment{violated buffers}
        \For{each $i \in B$}
            \State Select $j = \arg\min_{j' \in J(i)} v_{j'}$
            \State $z \gets x + |x_i| \cdot H_{:, j}$
            \If{all $z_i \geq -\epsilon$} \Comment{early stop check due to non-orthogonal reflection}
            \State \textbf{break}
            \Else
            \State $x \gets z$
            \EndIf
            \EndFor 
        \EndWhile
        \State $y \gets  \text{Solve } H y=(z - x)$
        \State \textbf{return} $z, y$
        \EndFunction
    \end{algorithmic}
\end{Subroutine}  
Our computational method seeks the neural network parameters ($w_1$, $w_2$) that minimize an approximation of the loss function defined in (\ref{eq:loss}). To do so, we first simulate discretized sample paths of the reference process $\tilde{Z}$ with the boundary pushing process $\tilde{Y}$ over time interval $[0, T]$; see Subroutine \ref{alg:euler}. This involves first sampling the paths of the underlying Brownian motion $X$. We then solve a discretized Skorokhod problem for each path of $X$ to obtain paths of $\tilde{Z}$ and $\tilde{Y}$. This is the purpose of Subroutine \ref{alg:skorokhod}. Thereafter, we compute a discretized version of the loss, summing over sampled paths to approximate the expectation and over discrete time steps to approximate the integral over $[0, T]$, and our method minimizes it using stochastic gradient descent; see Subroutine~\ref{alg:compute}.

\begin{Subroutine}[h]
\captionsetup{labelformat=empty}
    \caption{\textbf{Subroutine 3:} Computational Method}
    \label{alg:compute}
    \textbf{Input:} Number of steps $M$, batch size $B$, learning rate $r$, time horizon $T$, value vector $v_B = (v_1, \ldots, v_b)'$, discretization step size $h$ (assume $N = T/h$ is an integer), initial state $z$, and optimization solver ADAM
    
    \textbf{Output:} Neural network approximations $V_{w_1}$ (value function) and $G_{w_2}$ (gradient)
    \begin{algorithmic}[1]
        \State Initialize $V_{w_1}$ and $G_{w_2}$; set $z^{(i)}_0 = z$ for $i = 1, \ldots, B$

    \For{$k \gets 0$ to $M - 1$}
        \For{$i = 1$ to $B$}
\State \parbox[t]{\dimexpr\linewidth-\algorithmicindent}{
Simulate the Brownian increments $\{\tilde{Z}^{(i)}, \Delta \tilde{Y}^{(i)}, \delta^{(i)}\}$ via \textsc{Discretize}$(T, h, z_k^{(i)})$\\
\hspace*{1.5em}\Comment{See Subroutine \ref{alg:euler}}
} \EndFor
        \State Compute empirical loss
                \begin{align*}
                \hat{\ell}(w_1, w_2) = \frac{1}{B} \sum_{i=1}^B  &\Bigg( e^{-rT} V_{w_1}(\tilde{Z}^{(i)}(T)) - V_{w_1}(\tilde{Z}^{(i)}(0)) + \sum_{j=0}^{N-1} e^{-rhj} v_B \cdot \Delta \tilde{Y}^{(i)}(hj) \\
                & -  \sum_{j=0}^{N-1} e^{-rhj} G_{w_2}(\tilde{Z}^{(i)}(hj)) \cdot \delta_j^{(i)} + \sum_{j=0}^{N-1} e^{-rhj} F(\tilde{Z}^{(i)}(hj), G_{w_2}(\tilde{Z}^{(i)}(hj))) h \Bigg)^2 
                \end{align*}
        \State Compute gradients $\partial \hat{\ell} / \partial w_1$ and $\partial \hat{\ell} / \partial w_2$
        \State Update $w_1$, $w_2$ using ADAM
        \State Update $z_{k+1}^{(i)} \gets \tilde{Z}^{(i)}(T)$ for all $i$
    \EndFor
    
        \State \Return neural network approximations $V_{w_1}(\cdot)$ and $G_{w_2}(\cdot)$
    \end{algorithmic}
\end{Subroutine}
After obtaining the neural network parameters, our proposed policy is as follows: for $j = 1, \ldots, J$,
\begin{equation}\label{eq:theta_hat}
\hat{\theta}_j(z; w_2) =
\begin{cases}
  -\eta & \text{if } \hat{U}_j(z) \geq 0, \\
  \min\{\sqrt{n} x_j^*, \eta\} & \text{if } \hat{U}_j(z) < 0,
\end{cases}
\end{equation}
where $\hat{U}_j(z) = G_{\ell(j)}(z) + G_{r(j)}(z) + v_j$ for $z \in \mathbb{R}_+^I$.

\section{Proposed Matching Policy for the Pre-limit System}\label{sec:policy}

In this section, we interpret the policy (\ref{eq:theta_hat}) in the context of the pre-limit system by proposing a dynamic matching policy for that system. We begin with fixing a system parameter $n$ that will be used to scale various quantities of interest. Given the queue length vector $Q^n(t)$ in the $n$th system, the (scaled) queue length vector in the Brownian control problem is
\[
Z^n(t) = \frac{Q^n(t)}{\sqrt{n}}, \quad t \geq 0.
\]

In the $n$th system, consider the times at which either a job arrives to the system or a job waiting in the system abandons. We denote these times by $t_1, t_2, \ldots$ and refer to them as the decision epochs; and we let $t_0 = 0$ for notational convenience. Under our proposed policy, the system manager makes decisions only at these times. Because we assume all matches happen instantaneously, it is helpful to clarify the sequence of events at a decision epoch, say $t_k$. First, either an abandonment or arrival event occurs and the system state is updated. For notational convenience, we denote this updated (pre-decision) system state at time $t_k$ by $Z^n(t_k-)$ (and $Q^n(t_k-)$ for the unscaled version). Then the matching decisions are made and executed, and the system state is updated. We denote the post-decision state by $Z^n(t_k)$ (and $Q^n(t_k)$).
Recall that the system manager’s policy is described by the matching process $T^n(\cdot)$. Motivated by (\ref{eq:tb})-(\ref{eq:tnb}) and (\ref{eq:theta_hat}), we set $T^n(0) = 0$ and for $k \geq 1$, we let
\begin{equation}\label{eq:deltat}
\Delta T^n(t_k) = n x^* \Delta t_k - \sqrt{n} \int_{t_{k-1}}^{t_k} \theta(Z^n(s))\,ds - \sqrt{n} \Delta Y(t_k), 
\end{equation}
where $\Delta t_k = t_k - t_{k-1}$, $\Delta T^n(t_k) = T^n(t_k) - T^n(t_{k-1})$, and $\Delta Y(t_k) = Y(t_k) - Y(t_{k-1})$. Because $\Delta t_k$ is small, i.e., $O(1/ n)$, in the asymptotic regime we are focusing on, we approximate the integral term in (\ref{eq:deltat}) as follows:
\begin{equation}\label{eq:int_approx}
\int_{t_{k-1}}^{t_k} \theta(Z^n(s))\,ds \approx \theta(Z^n(t_k-)) \Delta t_k.
\end{equation}
To repeat, $Z^n(t_k-)$ denotes the (scaled) queue length vector immediately after the arrival or abandonment event at time $t_k$, but before any matches are made. Consequently, we have that
\begin{equation}\label{eq:deltaT}
\Delta T^n(t_k) \approx \left(n x^* - \sqrt{n}\, \theta(Z^n(t_k-))\right) \Delta t_k - \sqrt{n} \Delta Y(t_k).    
\end{equation}
Next, we define the \emph{potential number of matches} at time $t_k$, denoted by $\tau(t_k)$, as follows:
$$
\tau(t_k) = \left(n x^* - \sqrt{n} \, \theta(Z^n(t_k-))\right) \Delta t_k, \quad k \geq 1.
$$
Using (\ref{eq:theta_hat}), the proposed policy in terms of $\tau(\cdot)$ can be written as
\begin{align}\label{eq:tau}
\tau_j(t_k) = 
\left\{\begin{array}{ll}(nx_j^* + \sqrt{n}\eta)\Delta t_k, & \text{ if }\hat{U}_j(Z^n(t_k-))\geq 0, \\  (nx_j^* - \sqrt{n}\eta)^+\Delta t_k, & \text{ if } \hat{U}_j(Z^n(t_k-)) < 0,	\end{array}\right.
	\end{align}	 
where $a^+ = \max\{0, a\}$. Here, $\tau_j(t_k)$ denotes the intended number of matches at time $t_k$ via activity $j$ provided there are enough jobs in classes $\ell(j)$ and $r(j)$. Thus, the actual number of matches can be fewer. The purpose of the term $\sqrt{n} \Delta Y(t_k) $ in Equation (\ref{eq:deltaT}) is precisely to capture this. Before we address this issue, we turn to another implementation issue. The intended number $\tau_j(t_k)$ of matches via activity $j$ at time $t_k$ need not be an integer, whereas the actual matches must be integer valued. Defining $m_j(t_k)$ as the integer-valued planned number of matches via activity $j$ at time $t_k$, we set $$               
m_j(t_k) = \left\lceil \tau_j(t_k) - \varepsilon \right\rceil,
$$
where $\varepsilon \in [0,1]$ is a tuning parameter. In particular, setting $\varepsilon = 0, 0.5$, or $1$ corresponds to taking the ceiling, rounding, or taking the floor of the fractional argument $\tau_j(t_k)$. We let $m(t_k) = (m_j(t_k))$ be the vector of planned number of matches at time $t_k$ for $k \geq 1$.

Next, we turn to addressing the possibility of not having sufficiently many jobs to execute the desired number of matches. Although one can adopt the Skorokhod mapping approach to determine actual matching decisions, in principle, we proceed with a simpler approach. In doing so, we capture the effect of the pushing process $\Delta Y$ implicitly. More specifically, the proposed policy first performs the most valuable matches until either the intended number of matches is reached or one runs out of jobs to match using that activity. Formally, we let $a_1,\ldots, a_J$ be a permutation of $\{1, \ldots, J\}$ such that $\hat{U}_{a_1}(Z^n(t_k-)) \geq \cdots \geq \hat{U}_{a_J}(Z^n(t_k-))$. The proposed policy first focuses on activity $a_1$, i.e., the one with the highest shadow value, and performs
$$
\min\left\{ m_{a_1}(t_k), Q_{\ell(a_1)}^n(t_k), Q_{r(a_1)}^n(t_k) \right\}
$$
matches at decision epoch $t_k$. Once these matches are made, we update the state vector accordingly. Next, we consider activity $a_2$ and proceed similarly. Further details of the proposed matching policy are spelled out in Algorithm \ref{alg:static_review}.

\begin{MainAlgorithm}[h]
		\caption{Dynamic Matching Policy}
		\label{alg:static_review}
		\textbf{Input:} Queue length $ Q(t) \in \mathbb{R}^I $ at review time $t $; time increment $\Delta t$; optimal SPP solution $x^*$; scaling constant $n$; constants $\eta$, $\varepsilon$;  function $\hat{U}(\cdot): \mathbb{R}^I \to \mathbb{R}^J$ \\
		\textbf{Output:} A sequence of matching actions.
		\begin{algorithmic}[1]
			\State Initialize $ q \gets Q(t)$
             \State Compute $\hat{U} \gets (\hat{U}_1(q/\sqrt{n}), \ldots, \hat{U}_J(q/\sqrt{n}))$
		  \State Obtain index $(a_1,\ldots, a_J)$ s.t. $\hat{U}_{a_1} \geq \cdots \geq \hat{U}_{a_J}$
			\For{$ j \gets a_1\text{ to } a_J $}
			\If{$ \hat{U}_j \geq 0 $}
         \State $\tau_j = (nx_j^* + \sqrt{n}\eta)\Delta t$
          \Else \State $\tau_j = (n x_j^* - \sqrt{n} \eta)^+\, \Delta t$ \Comment{$\tau_j \equiv 0$ for $ j=b+1,\ldots, J$}
            \EndIf
            \State Compute $m_j = \lceil \tau_j-\varepsilon \rceil$
            \State Compute $d_j= \min(m_j, q_{l(j)}, q_{r(j)})$ \Comment{min over guided and available jobs}
            \If {$d_j \geq 1$}
			\State Perform $d_j$ matches on activity type $ j $\Comment{Execute matching}
            \State Update $q_{l(j)} \gets q_{l(j)} - d_j$, $q_{r(j)} \gets q_{r(j)} - d_j$
            \EndIf
			\EndFor
		\end{algorithmic}
\end{MainAlgorithm}
We also consider a variant of the proposed policy that we refer to as the dynamic matching policy with further updates. As done in the previous policy, at each decision point it considers activities sequentially in the order of their shadow values, i.e., in the order of $\hat{U}_j(Z^n(t_k-))$ at each decision epoch $t_k$. After performing the matches for the activity with the highest shadow value, the state is updated. The only difference between this policy and the previous one is that before considering which activity to focus on next among the remaining activities, their shadow values $\hat{U}_j$ are updated using the updated state vector. Details of this algorithm are provided in Algorithm \ref{alg:dynamic_review}.
\begin{MainAlgorithm}[H]
		\caption{Dynamic Matching Policy (with further updates)}
		\label{alg:dynamic_review}
		\textbf{Input:} Queue length $ Q(t) \in \mathbb{R}^I $ at review time $t $; time increment $\Delta t$; optimal SPP solution $x^*$; scaling constant $n$; constants $\eta$, $\varepsilon$;  function $\hat{U}(\cdot): \mathbb{R}^I \to \mathbb{R}^J $ \\
		\textbf{Output:} A sequence of matching actions.
		\begin{algorithmic}[1]
			\State Initialize $ q \gets Q(t) $
            \State Compute $\hat{U} \gets (\hat{U}_1(q/\sqrt{n}), \ldots, \hat{U}_J(q/\sqrt{n}))$
             \State Obtain index $(a_1,\ldots, a_J)$ s.t. $\hat{U}_{a_1} \geq \cdots \geq \hat{U}_{a_J}$
			\While{a matching activity is available}
			\For{$ j \gets a_1\text{ to }a_J $}
			\If{$ \hat{U}_{j} \geq 0 $}
                \State $\tau_j = (nx_{j}^* + \sqrt{n}\eta)\Delta t$
        \Else \State $\tau_j = (n x_j^* - \sqrt{n} \eta)^+\, \Delta t$ \Comment{$\tau_j \equiv 0$ for $ j=b+1,\ldots, J$}
            \EndIf
            \State Compute $m_j = \lceil \tau_j-\varepsilon \rceil$
            \State Compute $d_j= \min(m_j, q_{l(j)}, q_{r(j)})$ \Comment{min over guided and available jobs}
            \If {$d_{j} \geq 1$}
			\State Perform only one match on activity $ j $ \Comment{Execute matching}
            \State Update $q_{l(j)} \gets q_{l(j)} - 1$, $q_{r(j)} \gets q_{r(j)} - 1$
            \State Update $\hat{U} \gets (\hat{U}_1(q/\sqrt{n}), \ldots, \hat{U}_J(q/\sqrt{n}))$
			\State Update index $(a_1,\ldots, a_J)$ s.t. $\hat{U}_{a_1} \geq \cdots \geq \hat{U}_{a_J}$       
			\State \textbf{break} \Comment{Exit the inner for loop to restart the while loop}
            \EndIf
			\EndFor
			\EndWhile
		\end{algorithmic}
\end{MainAlgorithm}

\section{Benchmark Policies}\label{sec:benchmark}

To assess the performance of our proposed policy, we turn to benchmark policies drawn from the literature. More specifically, we consider the five benchmark policies described below. The first four make decisions at each decision epoch when the system state changes. The last benchmark, the static-priority policy, is a discrete-review policy, where the review interval length is a tunable hyperparameter. 

\noindent\textbf{Greedy policy.}
   The greedy policy executes all available matches immediately upon arrival. That is, under the greedy policy, when a job arrives, if there is a job that can be matched with it, then the match is executed at once and both jobs depart the system. If there are multiple possible matches, then the system manager chooses the one that generates the highest value. This policy doesn’t discriminate between basic and nonbasic activities. However, because nonbasic activities are inherently inefficient they should be used sparingly. The greedy policy may overuse them, leading to suboptimal performance; see, for example, the Zigzag A model in Section \ref{sec:results}. Thus, we also consider a version of the greedy policy that only uses the basic activities. We refer to that policy as the greedy-basic policy and describe it next.

\noindent\textbf{Greedy-basic policy.}
The greedy-basic policy executes only the matches associated with basic activities when they are available. However, the exclusion of nonbasic activities can be limiting because the use of certain nonbasic activities sparingly can improve system performance. For example, for the Zigzag C model of Section \ref{sec:results}, the greedy-basic policy underperforms because it never uses any nonbasic activities.

\noindent\textbf{First-come-first-serve (FCFS) policy.} Under the FCFS policy, each job is matched immediately upon arrival if there is a job it can be matched with. And, if there are multiple such jobs, then it is matched with the one that has been waiting the longest. Because FCFS policy considers only the order of arrival and ignores the value generated from the match, it can perform poorly. However, it can also outperform the greedy policies in some cases – particularly when the bipartite graph isn’t sufficiently flexible and the congestion costs are high. For example, see the Zigzag B model in Section~\ref{sec:results}.
         
\noindent\textbf{Longest-queue-first-serve (LQFS) policy.}
Under the LQFS policy, when a job arrives, if there is a job it can be matched with, then the match is executed immediately. If there are multiple possible matches corresponding to different buffers, the system manager chooses the job whose buffer is the longest. This policy can be viewed as a max-weight policy where the weights correspond to the number of jobs in the buffers, as done in \cite{VarBumMagWan2023}. Similar to FCFS, the LQFS policy does not take into account matching values. Nor does it make a distinction between basic versus nonbasic activities. Thus, one expects FCFS and LQFS to perform similarly.

\noindent\textbf{Static-priority policy.}
The static-priority policy is a discrete-review  policy that executes matches in a fixed priority order at pre-specified review epochs. At each review epoch, the policy first executes all available basic activities and then uses nonbasic activities to match any remaining jobs. To determine the priority order among the basic activities, we adopt the algorithm from \cite{AveDeLeviWard2024}, in which the static priority sets are obtained through a tree-pruning approach. For completeness, we restate the algorithm in Subroutine \ref{alg:priorityset} using our notation. For a concrete example to illustrate the policy, consider the model displayed in Figure \ref{fig:zz_c}. There are seven basic activities (denoted by solid lines connecting the buffers) ordered from top to bottom. At each review time, the highest priority is given to the outer most basic activities in the graph, i.e., the top activity (activity 1) and the bottom activity (activity 7). Any potential matches with these activities are made at once. Next, focusing attention on the rest of the graph, the system manager prioritizes activities 3 and 6. Lastly, she prioritizes the inner most activities, i.e., activities 4 and 5. After making all possible matches via basic activities, if there are jobs that can be matched by nonbasic activities, then the system manager executes them until no further matches can be made. The matches made via nonbasic activities can be executed in any order.

\begin{Subroutine}[h]
\captionsetup{labelformat=empty}
\caption{\textbf{Subroutine 4:} Static Priority Set Construction}
\label{alg:priorityset}
\textbf{Input:} Arrival rates $\lambda \in \mathbb{R}_+^{I}$; matching matrix $R$; optimal solution $x^*$\\
\textbf{Output:} A collection of priority sets $\{\mathcal{P}_0(x^*), \mathcal{P}_1(x^*), \ldots\}$
\begin{algorithmic}[1]
\State Initialize edge set $E_0 \gets \{(\ell(j),r(j)) \in \mathcal{E} : x^*_j > 0 \}$, $E \gets E_0$
\State Set capacities: $\mu \gets \lambda$, counter $h \gets 0$
\While{$E \neq \emptyset$}
    \State Initialize priority set: $\mathcal{P}_h(x^*) \gets \emptyset$
    \State Set consideration set $C \gets E$
    \While{$C \neq \emptyset$}
        \State Select any edge $j \in C$
        \If{$x^*_j = \mu_{\ell(j)}$ or $x^*_j = \mu_{r(j)}$}
            \State $\mathcal{P}_h(x^*) \gets \mathcal{P}_h(x^*) \cup \{j\}$
            \State Update: $\mu_{\ell(j)} \gets \mu_{\ell(j)} - x^*_{j}$, $\mu_{r(j)} \gets \mu_{r(j)} - x^*_{j}$
            \State Remove $(\ell(j),r(j))$ and neighbors from $C$: $C \gets C \setminus \{(\ell',r') : \ell'=\ell(j) \text{ or } r'=r(j) \}$
        \Else
            \State $C \gets C \setminus \{(\ell(j),r(j))\}$
        \EndIf
    \EndWhile
    \State Remove $\mathcal{P}_h(x^*)$ from $E$: $E \gets E \setminus \mathcal{P}_h(x^*)$
    \State $h \gets h + 1$
\EndWhile
\State Last priority set: $\mathcal{P}_h(x^*) \gets \{(\ell(j),r(j))\in \mathcal{E} : x^*_j = 0 \}$
\end{algorithmic}
\end{Subroutine}
It is worth noting that when the system only has basic activities, the static-priority policy resembles the greedy-basic policy provided the review period is sufficiently short. However, the two can still differ in terms of how they prioritize basic activities. The greedy-basic policy favors matchings that yield the highest immediate value, whereas the static-priority policy prioritizes matchings that incur the least opportunity cost, i.e., those with the smallest impact on the utilization of other basic activities, by starting with the top and the bottom activities and moving toward the center of the bipartite graph. Additionally, the static-priority policy allows the use of nonbasic activities to match any remaining jobs after all basic activities have been executed\footnote{The basic activities in our setting correspond to the activities characterized by the optimal extreme points in \cite{AveDeLeviWard2024}. They also consider another policy which is dominated by the static policy in their numerical experiments. Thus, we focus attention on the static priority policy.}.

\section{Test Problems}\label{sec:testproblems}
As mentioned in the introduction, we consider nine test problems. First, we consider three instances of the aforementioned X model; see Figure \ref{fig:x}. These instances differ in their abandonment rates. Second, we consider three 8-dimensional models that have $4\times4$ bipartite graphs. We refer to them as the Zigzag models. Third, we consider two test problems of dimension 24 that are constructed by concatenating the Zigzag models. Lastly, we consider a 120-dimensional test problem (with a $60\times60$ bipartite graph) that is constructed similarly. For the test problems, we assume $Z(0)=z=0$, and we set the scaling parameter $n=100$ and the discount rate $r=0.01$.
\subsection{X Models}
       	\begin{figure}[h]
		\centering
\includegraphics[width=0.4\textwidth]{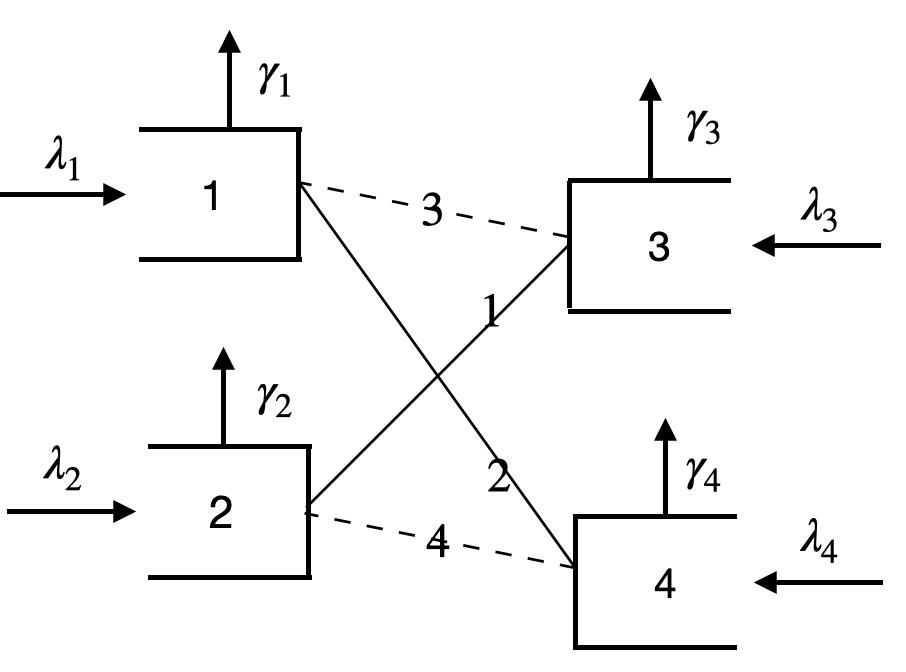}
		\caption{The X Model.}
		\label{fig:x}
	\end{figure}
The X model in Figure \ref{fig:x} has $2$ left classes and $2$ right classes.  Assume the model parameters are as follows: 
	\begin{align}
	\lambda= \begin{bmatrix} 0.5\\ 1\\ 1\\0.5\end{bmatrix}, \quad v = \begin{bmatrix} 4\\0.2\\2.0\\2.0\end{bmatrix}, \quad R = \begin{bmatrix}  0& 1&1 &0\\ 1& 0&0&1\\ 1&0&1&0\\0&1&0&1\end{bmatrix},\quad h= \begin{bmatrix} 1\\ 1\\ 1\\1\end{bmatrix},\quad a= 0,\quad \gamma = \alpha \begin{bmatrix} 1\\ 1\\ 1\\1\end{bmatrix},
	\end{align}
   where we consider three variants of the X model by setting $\alpha\in\{0.1, 0.01, 0.001\}$, which we refer to as the high, medium, and low abandonment levels, respectively. Solving the static planning problem (\ref{eq:spp}) yields the unique optimal solution $x^* = (1, 0.5, 0, 0)'$. 
	Accordingly, the basic activity set is $\left\{1, 2\right\}$, and the nonbasic set is $\left\{3, 4\right\}$. We then derive the proposed policy for the X model as described in Section \ref{sec:policy} and present its performance in Section \ref{sec:results}; also see Appendix \ref{app:NNarchi}, Table~\ref{table:parX} for the implementation details of our computational method.

\subsection{Zigzag Models}
    

\begin{figure}[htbp]
    \centering
    \begin{subfigure}[b]{0.3\textwidth}
        \centering
        \includegraphics[width=0.9\textwidth]{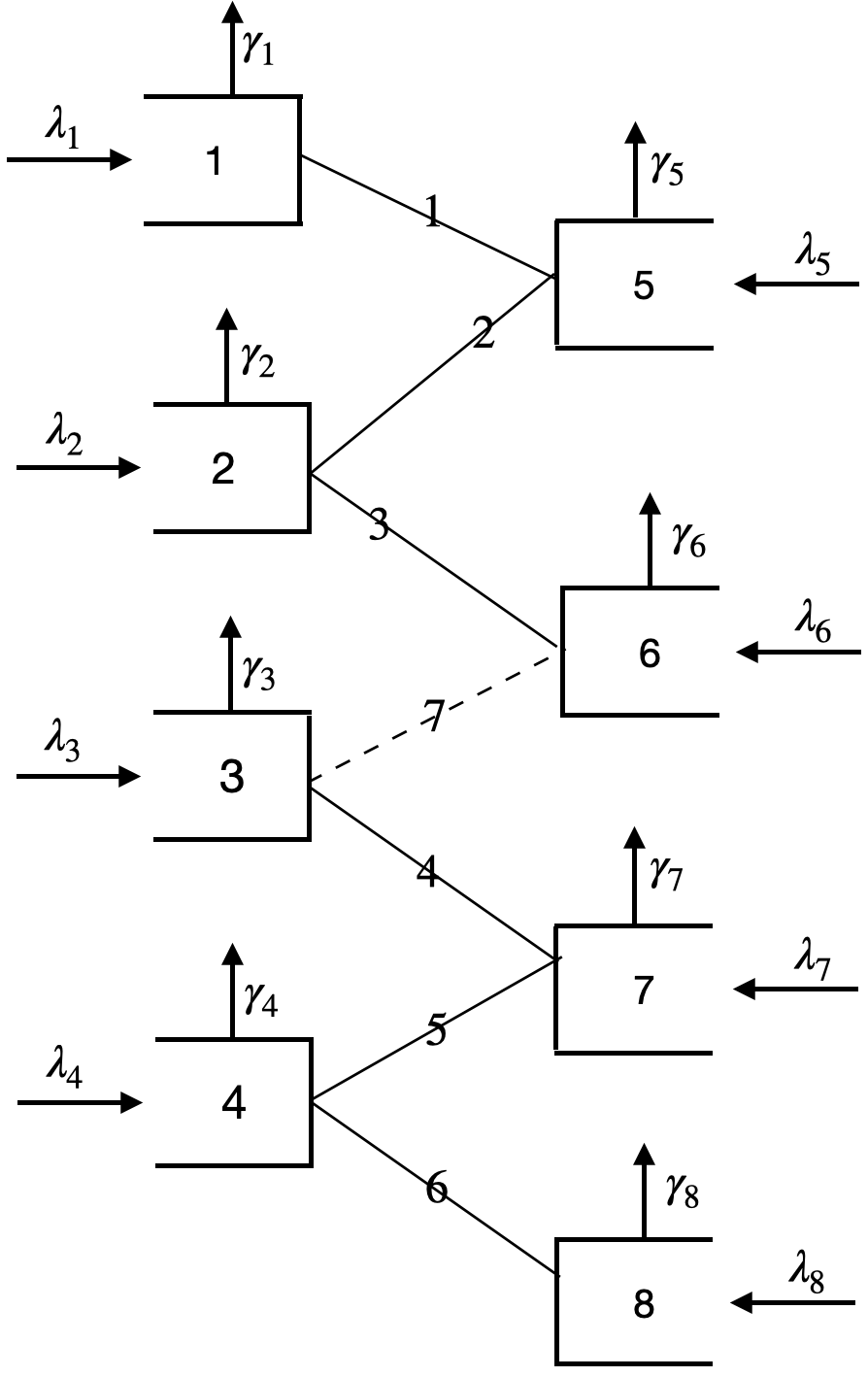}
        \subcaption{Zigzag A}
        \label{fig:zz_a}
    \end{subfigure}
    \hfill
    \begin{subfigure}[b]{0.3\textwidth}
        \centering
        \includegraphics[width=0.9\textwidth]{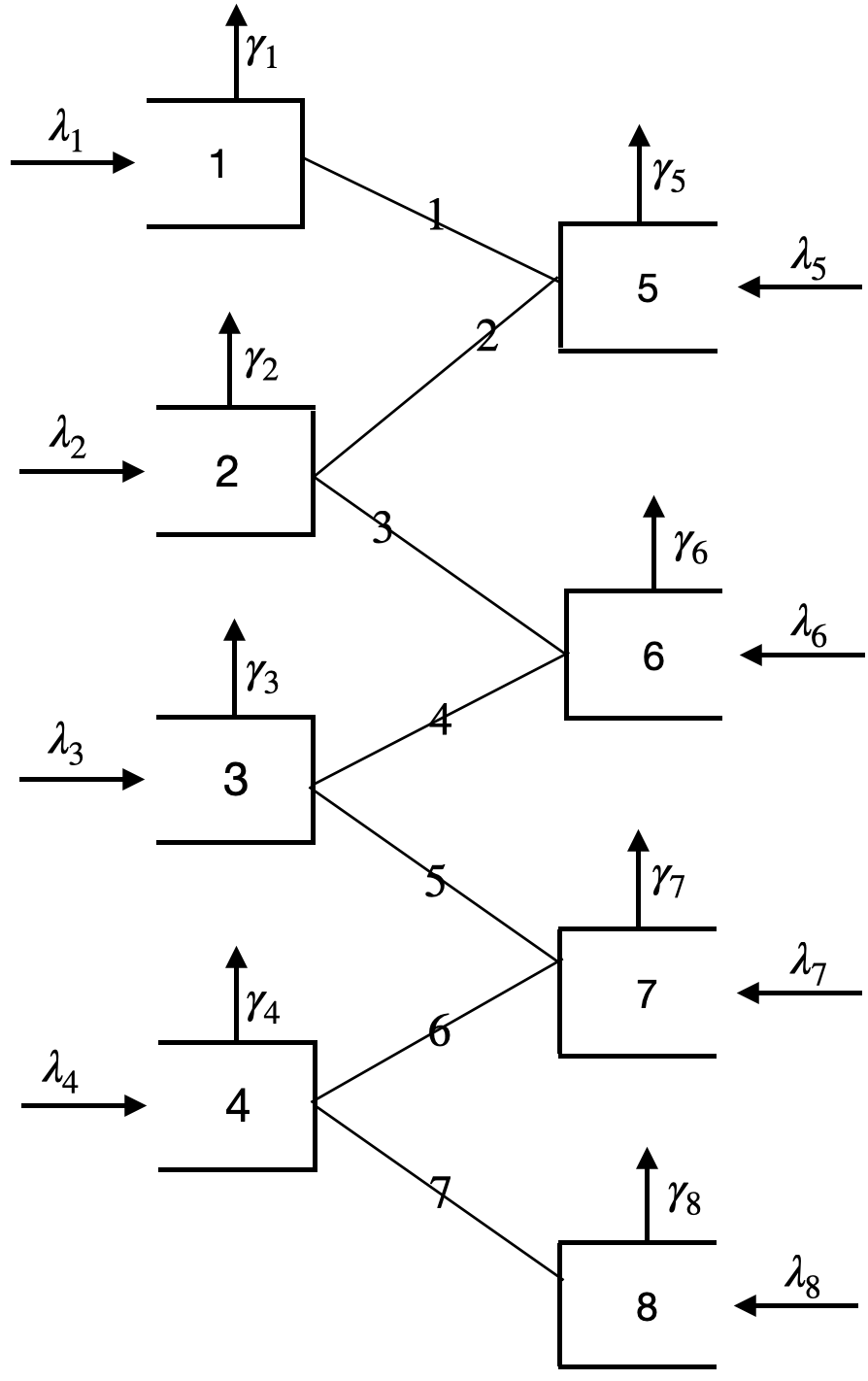}
        \subcaption{Zigzag B}
        \label{fig:zz_b}
    \end{subfigure}
    \hfill
    \begin{subfigure}[b]{0.3\textwidth}
        \centering
        \includegraphics[width=0.9\textwidth]{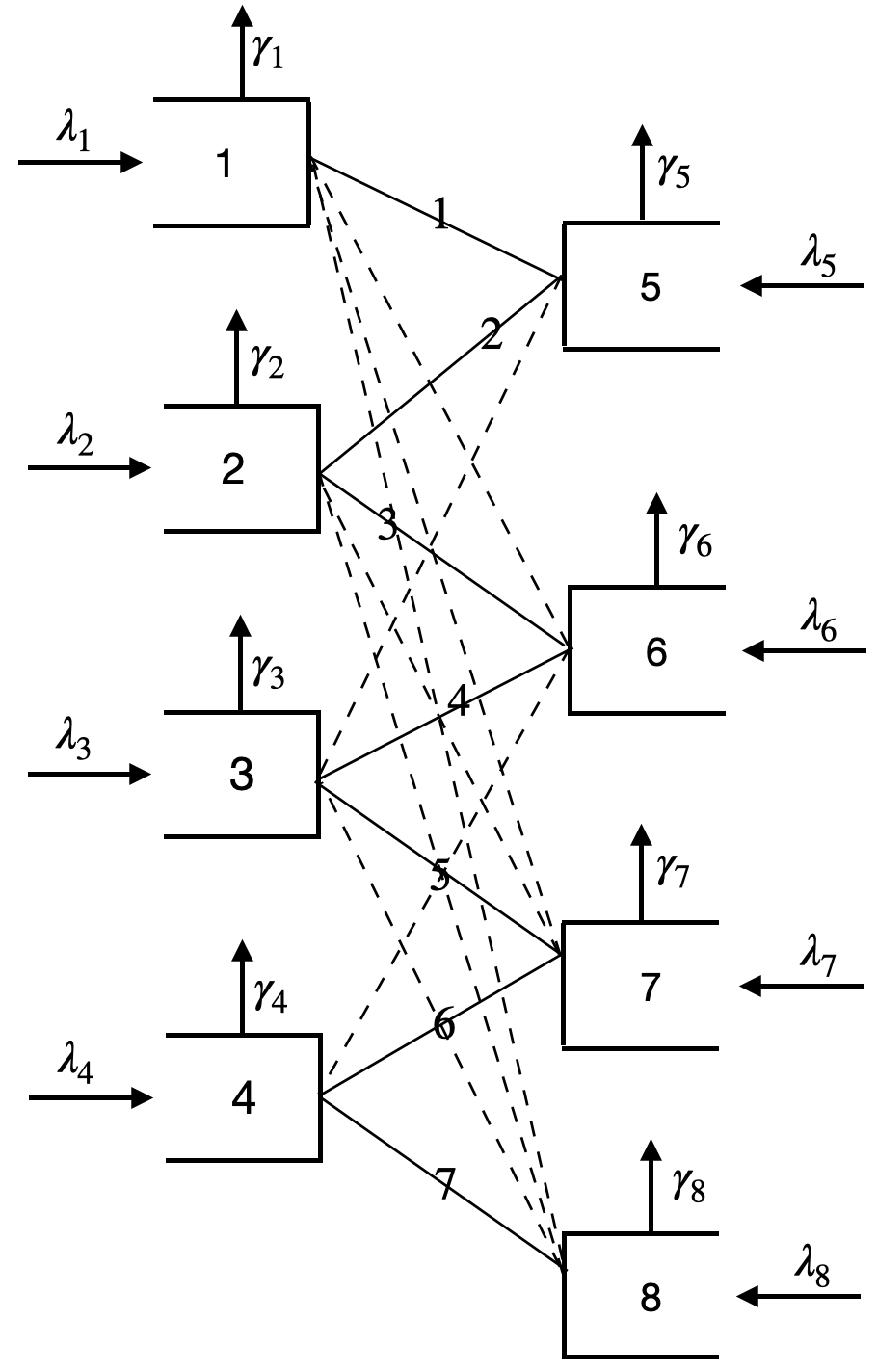}
        \subcaption{Zigzag C}
        \label{fig:zz_c}
    \end{subfigure}
    \caption{Zigzag models under different configurations.}
    \label{fig:zz_all}
\end{figure}

Our first test problem in this group involves a $4 \times 4$ bipartite network inspired by Figure 6 in \cite{KerAshGur2024}. It is reproduced in Figure~\ref{fig:zz_a} with slight modifications. The model parameters, along with the optimal solution $x^*$ to the static planning problem (\ref{eq:spp}), are summarized in Table~\ref{table:zz_par}. In particular, activities 1–6 are basic, whereas activity 7 is nonbasic. If one removes activity 7, the network decomposes into two separate sub-networks.
 
  To evaluate the performance of our model under different network structures, we design two additional variants: Zigzag B and Zigzag C. For Zigzag B, shown in Figure~\ref{fig:zz_b}, we construct a connected network composed of 7 basic activities. This is achieved by adjusting the parameters so that the solution $x^*$ (summarized in Table~\ref{table:zz_par}) includes only basic activities. Since Zigzag B consists solely of basic activities, one might expect that some simple policies could perform well—because no inefficient (nonbasic) activities are present. Due to its simplicity, this example is used to test whether our policy continues to offer significant improvements over the benchmark policies, and specifically whether the ordering or timing of basic matchings play a critical role.
     \begin{table}[H] 
    \centering
    \sisetup{
        separate-uncertainty = true,
        detect-all
    }
    {\renewcommand{\arraystretch}{1.1}
  {\small  \begin{tabular}{lccc}
        \toprule
        \multirow{2}{*}{} & \multicolumn{3}{c}{Value} \\
        \cmidrule(lr){2-4}
        & {Zigzag A} & {Zigzag B} & {Zigzag C} \\
        \midrule
        $\lambda$          & $(3, \,2,\, 1, \,4, \,4, \,1, \,2, \,3)'$  & $(3,\, 2,\, 2, \,4, \,4, \,2, \,2, \,3)'$   & $(3, \,2,\, 2, \,4, \,4, \,2, \,2, \,3)'$ \\
        $v$   & $(2, \,4, \,8, \,8, \,4, \,2, \,1 )'$  & $(2, \,4, \,8, \,9, \,8, \,4, \,2)'$  & $(2, \,4, \,8, \,9, \,8, \,4, \,2, \,2, \,4, \,2, \,6, \,2, \,2, \,4, \,2)'$\\
        $x^*$ & $(3, \,1, \,1,\, 1, \,1, \,3, \,0)'$ & $(3, \,1, \,1, \,1, \,1, \,1, \,3)'$   & $(3, \,1, \,1,\, \,1, \,1, \,1, \,3, \,0, \,0, \,0,\,0, \,0, \,0, \,0,\, 0)'$ \\
\midrule
        $R$   &$\begin{bmatrix}
1 & 0 & 0 & 0 & 0 & 0 & 0 \\
0 & 1 & 1 & 0 & 0 & 0 & 0 \\
0 & 0 & 0 & 1 & 0 & 0 & 1 \\
0 & 0 & 0 & 0 & 1 & 1 & 0 \\
1 & 1 & 0 & 0 & 0 & 0 & 0 \\
0 & 0 & 1 & 0 & 0 & 0 & 1 \\
0 & 0 & 0 & 1 & 1 & 0 & 0 \\
0 & 0 & 0 & 0 & 0 & 1 & 0
\end{bmatrix}$ & 
$\begin{bmatrix}
1 & 0 & 0 & 0 & 0 & 0 & 0 \\
0 & 1 & 1 & 0 & 0 & 0 & 0 \\
0 & 0 & 0 & 1 & 1 & 0 & 0 \\
0 & 0 & 0 & 0 & 0 & 1 & 1 \\
1 & 1 & 0 & 0 & 0 & 0 & 0 \\
0 & 0 & 1 & 1 & 0 & 0 & 0 \\
0 & 0 & 0 & 0 & 1 & 1 & 0 \\
0 & 0 & 0 & 0 & 0 & 0 & 1
\end{bmatrix}$
    & $\left[\begin{array}{@{}ccccccccccccccc@{}}
1 & 0 & 0 & 0 & 0 & 0 & 0 & 1 & 1 & 1 & 0 & 0 & 0 & 0 & 0 \\
0 & 1 & 1 & 0 & 0 & 0 & 0 & 0 & 0 & 0 & 1 & 1 & 0 & 0 & 0 \\
0 & 0 & 0 & 1 & 1 & 0 & 0 & 0 & 0 & 0 & 0 & 0 & 1 & 1 & 0 \\
0 & 0 & 0 & 0 & 0 & 1 & 1 & 0 & 0 & 0 & 0 & 0 & 0 & 0 & 1 \\
1 & 1 & 0 & 0 & 0 & 0 & 0 & 0 & 0 & 0 & 0 & 0 & 1 & 0 & 0 \\
0 & 0 & 1 & 1 & 0 & 0 & 0 & 1 & 0 & 0 & 0 & 0 & 0 & 0 & 1 \\
0 & 0 & 0 & 0 & 1 & 1 & 0 & 0 & 1 & 0 & 1 & 0 & 0 & 0 & 0 \\
0 & 0 & 0 & 0 & 0 & 0 & 1 & 0 & 0 & 1 & 0 & 1 & 0 & 1 & 0
\end{array}\right]$ \\
        \midrule
         $h$ & \multicolumn{3}{c}{$(11.2, \,12.9, \,9.6, \,3.5,\, 0, \,19.5, \,11.4, \,1.3)'$} \\
        $a$ & \multicolumn{3}{c}{$(8.0, \,0.9, \,15.3, \,19.6, \,2.0, \,1.0, \,0.6, \,10.2)'$} \\
        $\gamma$ & \multicolumn{3}{c}{$(0.094, \,0.146,\, 0.061,\, 0.177, \,0.082, \,0.143, \,0.053,\, 0.049)'$} \\
        \bottomrule
    \end{tabular}}}
    \caption{Parameters of the Zigzag models.}
    \label{table:zz_par}
\end{table}       
   Zigzag C, shown in Figure~\ref{fig:zz_c}, is almost fully connected model. In addition to the 7 basic activities used in Zigzag B, Zigzag C includes 8 nonbasic activities, making the matching problem more challenging. While nonbasic activities may help reduce holding and abandonment costs in some instances, their excessive use can limit the use of more efficient basic activities. The model parameters and the optimal solution $x^*$ are again summarized in Table~\ref{table:zz_par}.
   
  The holding costs, abandonment costs, and abandonment rates were not specified in \cite{KerAshGur2024}; therefore, we use randomly generated values. Specifically, both holding and abandonment cost parameters are drawn uniformly from $[0, 20]$, and abandonment rates are drawn uniformly from $[0, 0.2]$. All Zigzag models share the same holding costs, abandonment rates, and abandonment cost parameters, which are summarized in Table~\ref{table:zz_par}. Furthermore, model training parameters are summarized in Appendix~\ref{app:NNarchi}, Table~\ref{table:parZZ}.   

\subsection{High Dimensional Models}
To show the effectiveness of our computational method in higher dimensions, we consider two 24-dimensional models and a 120-dimensional model. The first 24-dimensional model is constructed by concatenating the three Zigzag models—thus incorporating diverse network features—and adding 3 basic activities. More specifically, the bipartite network of this model is connected through basic activities. For the second 24-dimensional model, we concatenate copies of the Zigzag C model using three nonbasic activities.

The 120-dimensional model is formed by concatenating five copies of the second 24-dimensional model introducing 4 basic activities to connect them. The resulting system has 120 buffers, 109 basic activities and 135 nonbasic activities. All test problems are designed to satisfy Assumption \ref{as:bf}. Detailed model parameters, including the arrival rates, the matching matrix R, randomly generated abandonment rates, holding costs, and abandonment costs, are provided in Appendix \ref{app:hdtest}.  We summarize the model training parameters in  Appendix \ref{app:NNarchi}, Table~\ref{table:parHD}.

\section{Computational Results}\label{sec:results}
\begin{figure}[H]
		\centering
\includegraphics[width=0.7\textwidth]{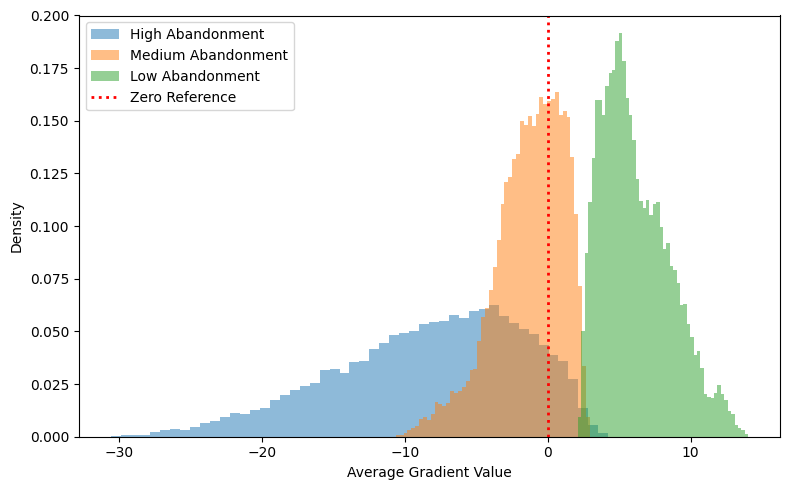}
		\caption{Distribution of Average Gradients Across Abandonment Levels}
		\label{fig:histGrad}
	\end{figure}

In this section, we compare our two policies with the benchmark policies using Monte Carlo simulation (using common random numbers across different policies). Throughout, we use the centered discounted value as the evaluation criterion. Thus, a smaller value indicates a better performance. Tables \ref{table:perf_x}-\ref{table:perf_24_120} report 95\% confidence intervals. In these tables, "our policy" refers to Algorithm~\ref{alg:static_review} whereas "our policy (further updates)" refers to Algorithm~\ref{alg:dynamic_review}. Also, for the static-priority policy, the review period length $l$ is determined via a one-dimensional grid search; see Table \ref{table:review_period} in Appendix \ref{app:NNarchi} for the resulting values.
\paragraph{X model.}
Table \ref{table:perf_x} provides the performance comparison of our policies with the benchmark policies. Our policies significantly outperform the other benchmark policies.
  \begin{table}[h] 
    \centering
    \sisetup{
        separate-uncertainty = true,
        detect-all
    }
    \begin{tabular}{l
                    S
                    S
                    S}
        \toprule
        \multirow{2}{*}{Policy} & \multicolumn{3}{c}{\textbf{Centered discounted value} (in hundreds)} \\
        \cmidrule(lr){2-4}
& {\textbf{\shortstack{High \\ Abandonment}}} 
& {\textbf{\shortstack{Medium \\ Abandonment}}} 
& {\textbf{\shortstack{Low \\ Abandonment}}} \\
        \midrule
        Greedy          & \num{3.63 \pm 0.10}   & \num{7.51\pm 0.43}   & \num{10.55 \pm 0.81} \\
        Greedy-basic    & \num{4.68 \pm 0.10}   & \num{10.34 \pm 0.44}  & \num{14.79 \pm 0.83} \\
        FCFS            & \num{3.99 \pm 0.10}  & \num{7.70 \pm 0.44}   & \num{10.70 \pm 0.83} \\
        LQFS            & \num{3.90 \pm 0.10}  & \num{7.67 \pm 0.44}   & \num{10.67 \pm 0.83} \\
        Static-priority & \num{3.47 \pm 0.10}  & \num{7.16 \pm 0.43}   & \num{10.17 \pm 0.83} \\
        \midrule
        Our policy        & \num{3.34 \pm 0.10}   & \num{7.14 \pm 0.43}  & \num{10.12 \pm 0.81} \\
        Our policy (further updates)       & \num{3.31 \pm 0.10}   & \num{7.07 \pm 0.43}  & \num{10.08 \pm 0.81} \\
        \bottomrule
    \end{tabular}
    \caption{Performance comparison for the X models under different abandonment rates based on 100 simulation replications.}
    \label{table:perf_x}
\end{table}
Because the X model is relatively simple, one expects some benchmark policies to do well as we observe in Table \ref{table:perf_x}. To further explore this, we consider the number of times each activity is undertaken during the first thousand time units. For brevity, we focus on the high abandonment case (the other cases are similar). Recall that activities 1 and 2 are basic whereas activities 3 and 4 are nonbasic (see Figure~\ref{fig:x}). Our policy makes a small but non-negligible number of matches using nonbasic activities. One expects from Equations (\ref{eq:tb})-(\ref{eq:tnb}) that the basic activity rates are $O(n)$ whereas nonbasic activity rates are $O(\sqrt{n})$. Thus, we expect their ratio to be of order $\sqrt{n} = 10$, cf. Table \ref{table:AveActivity}. This is also true for the static-priority policy since it strives to match using the basic activities. Its performance is closest to that of our proposed policies. 
\begin{table}[h] 
    \centering
    \begin{tabular}{lccccc}
        \toprule
        \multirow{2}{*}{Policy} & \multicolumn{5}{c}{Average number of matches} \\
        \cmidrule(lr){2-6}
        & $1$ & $2$ & & $3$ & $4$ \\
        \midrule
        Greedy          & 67217 & 15880 & & 32681 & 32709 \\
        Greedy-basic    & 98719 & 49142 & & 0     & 0 \\
        FCFS            & 65968 & 16526 & & 32979 & 33012 \\
        LQFS            & 66268 & 16341 & & 32934 & 32943 \\
        Static-priority & 95772 & 46161 & & 3259  & 3285 \\
        \midrule
        Our policy       & 95011 & 44534 & & 4453  & 4439 \\
        Our policy (further updates)       & 94526 & 44078 & & 4933  & 4919 \\
        \bottomrule
    \end{tabular}
    \caption{A comparison of activity usage levels for the X model based on the first 1000 time units of the simulation.}
    \label{table:AveActivity}
    
\end{table}
\paragraph{Abandonment rates, market thickness, and their effect on shadow prices $\mathbf{\nabla V(\cdot)}$ and $\mathbf{G(\cdot, \boldsymbol{w}_2)}$.} 
To better understand how our policy performs under different levels of abandonment, we examine the gradients of the value function $\nabla V(\cdot)$, approximated by the trained neural network $G(\cdot; w_2)$. These gradients capture the sensitivity of the value function to changes in the system state, indicating how the value evolves as the number of jobs in the system increases. So they can be interpreted as the shadow prices in a certain sense. We compute the gradients with respect to the four queue lengths and plot the histogram of their average. That is, we compute the realized gradient estimates along the sample paths of the state trajectories observed in our simulation study under the proposed policy. In Figure~\ref{fig:histGrad}, we present these histograms for the three variants of the X model, corresponding to high, medium, and low abandonment levels. As shown in Figure \ref{fig:histGrad}, higher abandonment levels lead to more negative average gradients under our policy.

Recall that we set the abandonment cost to zero for the instances of X model we consider. Thus, there is no direct abandonment cost. The core tradeoff is between reducing the holding cost and increasing the number of matches, thereby increasing the value generated by those matches. When abandonment is low, jobs tend to accumulate, forming a thick market that facilitates successful matching but at the expense of high holding costs. In this regime, positive gradients reflect the fact that lowering the number of jobs in the system improves the system performance, because the holding costs tend to dominate the cost. In contrast, when the abandonment rates are high, few jobs remain in the system, resulting in a thin market. Matching opportunities become  more scarce due to frequent abandonments. Holding costs no longer dominate the system performance. In this case, the negative gradient reflects the fact that increasing the number of jobs in the system facilitates more matches and creates value, improving the system performance.

In essence, our proposed policy adapts to current market thickness dynamically. On the one hand, when it is deemed thick, the gradient is positive and the system manager strives to execute the feasible matches aggressively. On the other hand, when it is deemed thin, she may hold off on some feasible matches to accumulate more jobs in order to facilitate higher-value matches, cf. Equations (\ref{eq:theta*}) and (\ref{eq:tau}).

\paragraph{Zigzag models.}   
Table \ref{table:perf_zz} shows that our proposed policy significantly outperforms the benchmark policies for the Zigzag models with the following nuances. 
\begin{table}[H] 
    \centering
    \sisetup{
        separate-uncertainty = true,
        table-format=5.2(2),
        detect-all
    }
    \begin{tabular}{lSSS}
        \toprule
        \multirow{2}{*}{Policy} & \multicolumn{3}{c}{\textbf{Centered discounted value} (in thousands)} \\
        \cmidrule(lr){2-4}
        & {\textbf{Zigzag A}} & {\textbf{Zigzag B}} & {\textbf{Zigzag C}} \\
        \midrule
        Greedy & \num{39.74 \pm 0.42}~~  & \num{19.42\pm 0.47} &\num{9.91 \pm 0.20} \\
        Greedy-basic    & \num{5.73 \pm 0.27}~~    & \num{19.42 \pm 0.47} &\num{19.42 \pm 0.47} \\
        FCFS            & \num{9.98 \pm 0.23}~~    & \num{6.76 \pm 0.22} &\num{17.91 \pm 0.22} \\
        LQFS            & \num{11.24 \pm 0.25}~~   & \num{7.96 \pm 0.25} &\num{19.21 \pm 0.25} \\
        Static-priority & \num{14.41 \pm 0.39} ~~  & \num{12.60 \pm 0.39} &\num{14.75 \pm 0.39} \\
        \midrule
        Our policy        & \num{5.25 \pm 0.21}~~    & \num{4.15 \pm 0.21}  &\num{5.13 \pm 0.34} \\
        Our policy (further updates)       & \num{5.21 \pm 0.21} ~~   & \num{3.87\pm 0.20} &\num{4.90 \pm 0.34} \\
        \bottomrule
    \end{tabular}
    \caption{Performance comparison for Zigzag models based on 100 simulation replications.}
    \label{table:perf_zz}
\end{table}
First, for the Zigzag A model, the greedy-basic policy significantly outperforms the other benchmarks, and its performance is within 10\% of our proposed policies. As shown in Figure \ref{fig:zz_a}, if one deletes the only nonbasic activity (activity 7), the network decomposes into two simple sub-networks consisting of basic activities. Given that our proposed policies use this nonbasic activity only sparingly by design, and the greedy-basic policy never uses it as well as given the simplicity of the two sub-networks, we expect them to perform somewhat similarly. Further exploration of activity usage levels underscores this similarity; see Table \ref{table:AveActivity_zz}.

\begin{table}[h] 
    \centering
    \begin{tabular}{lcccccccc}
        \toprule
        \multirow{2}{*}{Policy} & \multicolumn{8}{c}{Average number of matches} \\
        \cmidrule(lr){2-9}
        & $1$ & $2$ &$3$ & $4$& $5$ & $6$ & & $7$ \\
        \midrule
        Greedy-basic  &   297346& 100085& 99894& 99921& 100050& 297399 &&0\\
        Our policy & 299593& 97606& 99911& 
99919& 100019& 297420&& 3 \\
        \bottomrule
    \end{tabular}
    \caption{A comparison of activity usage levels for Zigzag A model based on the first 1000 time units of the simulation.}
    \label{table:AveActivity_zz}
\end{table}
Second, in the Zigzag B model, the system is fully connected through basic activities, and no nonbasic activity is available. This means every available activity is efficient and should be utilized. In this case, reducing congestion costs becomes the main concern. As shown in Table~\ref{table:perf_zz}, the FCFS and LQFS policies perform best among the benchmarks. However, those policies ignore the differences in the value rates of different activities whereas our policies take them into account. Thus, our policies significantly outperform them.

Third, in the Zigzag C model, the system is more complex, with many nonbasic activities available. A trade-off emerges: using less efficient activities to reduce congestion costs versus using more efficient ones to increase matching value. In this example, the greedy policy performs best among the benchmarks, although this wasn't obvious to us apriori.

Despite the varying structural characteristics and which type of policy they favor, our policy consistently adapts to each scenario and significantly outperforms all benchmarks.

\paragraph{24- and 120-dimensional models.}
In the high-dimensional test problems, our policy consistently outperforms the benchmarks by significant margins. For the 24 dimensional models, where we simply concatenate the Zigzag models, the static-priority
policy performs the best among the benchmarks. Notably, as mentioned earlier, the static-priority policy performs worse than the other benchmarks when applied to any individual Zigzag model. In contrast, the best benchmark policy for the 120 dimensional example is the greedy-basic policy although it is significantly outperformed by our proposed policies.



\begin{table}[h] 
    \centering
    \begin{tabular}{lSSS}
        \toprule
        \multirow{2}{*}{Policy} & \multicolumn{3}{c}{\textbf{Centered discounted value} (in thousands)} \\
        \cmidrule(lr){2-4}
        & {\textbf{24-D (I)}} & \textbf{24-D (II)} & \textbf{120-D} \\
        \midrule
        Greedy            & \num{17.24 \pm 0.22}  & \num{55.75\pm 1.17} & \num{203.80 \pm 0.64} \\
        Greedy-basic      & \num{27.29 \pm 0.72}  & \num{56.09\pm 3.41} &\num{71.75\pm 0.64} \\
        FCFS            & \num{19.12 \pm 0.24} &\num{69.27 \pm 1.00} & \num{286.80 \pm 0.27} \\
        LQFS             & \num{19.42 \pm 0.26}  & \num{68.70\pm 1.00} &\num{285.51\pm 0.25} \\
        Static-priority    & \num{8.97 \pm 0.30}   &\num{31.72 \pm 1.33}& \num{92.81 \pm 0.99} \\
        \midrule
        Our policy   & \num{7.12 \pm 0.20}                       & \num{20.66 \pm 1.28} & \num{49.19\pm 0.52} \\
        Our policy (further updates)      & \num{6.85\pm 0.20}                         &\num{20.67\pm 1.28}& \num{49.19\pm 0.52} \\
        \bottomrule
    \end{tabular}
    \caption{Performance comparison for 24-D and 120-D models based on 30 simulation replications.}
    \label{table:perf_24_120}
\end{table}

    \section{Concluding Remarks} \label{sec:conclusion}
We design a dynamic control policy tailored for large-flow, high-dimensional matching systems. We approximate the control problem by deriving a Brownian control problem, a drift-rate control problem of a reflecting Brownian motion. We solve this problem by adopting the deep BSDE approach that relies heavily on deep neural network technology. Using this solution, we propose a matching policy that can be directly implemented in the original system.

Overall, our dynamic matching policy consistently outperforms the benchmark policies. Most benchmark policies make decisions solely in a myopic fashion\footnote{Although the static-priority policy can be slightly tuned by delaying the review time to accumulate more information, it still relies entirely on the system state at the time of review.}. For example, the next arriving job—potentially arriving shortly—might present a more valuable matching opportunity if the current job were held just a bit longer. On the other hand, it could be beneficial to perform a seemingly suboptimal activity immediately if doing so contributes to improved long-term system performance. Such decision-making requires knowledge of the value function. As demonstrated across various test problems in Section \ref{sec:results}, our policies perform better because we incorporate gradient of the value function, that is, we account for the impact of current decisions on the future evolution and performance of the system.


For simplicity, Assumption \ref{as:bf} requires that the system is perfectly balanced. This can be extended to allow second-order imbalances as is usually done in the heavy-traffic literature, see for example, Section 6 of \cite{Harr2003}. Also, our computational method can be applied to real-world matching problems, e.g., for making dynamic dispatch decisions in ride-hailing platforms. These extensions and applications are left for future work. 

\begin{singlespace}
\bibliographystyle{ims}
\bibliography{reference}
\end{singlespace}

	\appendix
\renewcommand{\theequation}{\Alph{section}.\arabic{equation}}
\setcounter{equation}{0}
	\begin{appendix}
		\addtocontents{toc}{\protect\setcounter{tocdepth}{1}}
		\makeatletter
		\addtocontents{toc}{%
			\begingroup
			\let\protect\l@chapter\protect\l@section
			\let\protect\l@section\protect\l@subsection
		}

    \section{A formal derivation of the Brownian Control Problem}\label{app:deriveBCP}
    This appendix provides the derivation of the Brownian control problem (BCP). For the $n$th system, denote by $E^n(t)$ the vector of arrival processes with rate $\lambda^n$ and by $A^n(t)$ the abandonment process. The evolution of the system state in the $n$th system is
    $$Q^n(t) = Q^n(0) + E^n(t) - R T^n(t) - A^n(t), \quad t \geq 0,$$
    where $T^n(t)$ satisfies (\ref{eq:tb})-(\ref{eq:tnb}) and $A^n(t) = N\left(\int_0^t \Gamma Q^n(s) d s\right)$. Then, by definition (\ref{eq:zn}), the scaled queue length process can be written as
    \begin{align}\label{eq:zn1}
    Z^n(t) = Z^n(0) +\frac{1}{\sqrt{n}}  E^n(t) - \frac{1}{\sqrt{n}} RT^n(t) - \frac{1}{\sqrt{n}} A^n(t).
    \end{align}
    By strong approximations for Poisson processes (\cite{HorCso1993}), for $i=1, \ldots, I$, we have the following expansions:
    $$\begin{aligned}
    &E_i^n(t) = \lambda_i nt + \sqrt{\lambda_i n}B_i(t) + o(\sqrt{n}),\\
    &A_i^n(t) = \int_0^t \gamma_i Q^n_i(s) d s +o(\sqrt{n}), 
    \end{aligned}$$
    where $\{B_i(t), t\geq 0\}$ for $i = 1,\ldots, I$ are independent standard Brownian motions. Substituting these expansions and (\ref{eq:tb})-(\ref{eq:tnb}) into (\ref{eq:zn1}), we have
    \begin{equation}\label{eq:zn2}
    \begin{aligned}
    Z^n(t)= &Z^n(0) + X(t) + \lambda \sqrt{n} t - Rx^* \sqrt{n}t + \int_{0}^tR\,\theta(
    Z^n(s))ds + H Y_B(t)- \int_0^t \Gamma Z^n(s) d s + o(1)\\
    =&Z^n(0) + X(t) + \int_{0}^tR\,\theta(
    Z^n(s))ds + H Y_B(t)- \int_0^t \Gamma Z^n(s) d s + o(1),
    \end{aligned}
    \end{equation}
    where we applied $Rx^* = \lambda$ by Assumption \ref{as:bf}, and $X(t)$ is an $I$-dimensional driftless Brownian motion with a diagonal covariance matrix $\Lambda = diag(\lambda_1,\ldots , \lambda_I )$. Passing to the limit formally in (\ref{eq:zn2}) as $n\rightarrow\infty$, we obtain the limiting queue length process (\ref{eq:z}) in the approximating Brownian control problem. 
    
  Similarly, for the scaled cost process (\ref{eq:xi_n}), we substitute $\Pi^n(t)$, the expansion of $A^n(t)$, and (\ref{eq:tb})–(\ref{eq:tnb}) into (\ref{eq:xi_n}) to obtain:
    \begin{align}\label{eq:xi1}
    \xi^n(t) &= \frac{1}{\sqrt{n}} \left\{(v \cdot x^*) n t - v\cdot T^n + \int_0^t h\cdot Q^n(s) d s + \int_0^t a\cdot \Gamma Q^n(s) ds \right\} + o(1)\\
    &= \int_{0}^tv\cdot \theta(
    Z^n(s))ds + v_B\cdot Y_B(t) + \int_0^t(h+a\Gamma)\cdot Z^n(s) ds+ o(1). 
    \end{align}
     Then passing to the limit formally as $n\rightarrow \infty$, we obtain the limiting cost process (\ref{eq:xi}) in the approximating Brownian control problem.
    
    \section{Proof of Proposition \ref{equivSDE}}\label{app:equivSDE}
\begin{proof}
	Applying It\^o's formula on $e^{-rt}V(\tilde{Z}(t))$, one has
	\begin{equation}\label{equiv_0}
	\begin{aligned}
	e^{-r T} V(\tilde{Z}(T))-V(\tilde{Z}(0)) & =\int_0^T e^{-r t} \nabla V(\tilde{Z}(t)) \cdot \mathrm{d} W(t) +\int_0^T e^{-r t} \mathcal{D}V(\tilde{Z}(t)) \cdot \mathrm{d} \tilde{Y}_{B}(t)\\
	& +\int_0^T e^{-r t} \left[\mathcal{L}V(\tilde{Z}(t)) + (R\tilde{\theta}-\Gamma\tilde{Z}(t))\cdot\nabla V(\tilde{Z}(t))-r V(\tilde{Z}(t)) \right] \mathrm{d} t.
	\end{aligned}
	\end{equation}
	Using the complementarity condition (\ref{eq:complem}) and the boundary condition (\ref{eq:hjbb}), one has
	\begin{align}\label{equiv_1}
	\int_0^T e^{-r t} \mathcal{D}V(\tilde{Z}(t)) \cdot \mathrm{d} \tilde{Y}_{B}(t) = -\int_0^T e^{-r t} v_B \cdot \mathrm{d} \tilde{Y}_{B}(t).
	\end{align}
	Substituting $\tilde{Z}(t)$ into the HJB equation yields
	\begin{align}
	\mathcal{L}V(\tilde{Z}(t)) + \min _{\theta \in \Theta}\{(R'\nabla V(\tilde{Z}(t)) + v)\cdot\theta \} +\big(c-\Gamma\nabla V(\tilde{Z}(t)) \big)\cdot \tilde{Z}(t)=r V(\tilde{Z}(t)).
	\end{align}
	Multiplying both sides by $e^{-rt}$ and integrating over $[0,T]$ yields
	\begin{align}\label{equiv_2}
	&\int_0^T e^{-r t} \left[\mathcal{L}V(\tilde{Z}(t))-r V(\tilde{Z}(t)) \right] \mathrm{d} t \\
    &= -\int_0^T e^{-r t}\left[\min _{\theta \in \Theta}\{(R'\nabla V(\tilde{Z}(t)) + v)\cdot\theta \} +\big(c-\Gamma\nabla V(\tilde{Z}(t)) \big)\cdot \tilde{Z}(t)\right] \mathrm{d} t
	\end{align}
	Plugging (\ref{equiv_1}) and (\ref{equiv_2}) in (\ref{equiv_0}), we have
	$$
	\begin{aligned}
	e^{-r T} V(\tilde{Z}(T))-V(\tilde{Z}(0))  =&\int_0^T e^{-r t} \nabla V(\tilde{Z}(t)) \cdot \mathrm{d} W(t) -\int_0^T e^{-r t} v_B \cdot \mathrm{d} \tilde{Y}_{B}(t)\\
	&+\int_0^T e^{-r t}(R\tilde{\theta}-\Gamma\tilde{Z}(t)) \cdot \nabla V(\tilde{Z}(t))\mathrm{d} t\\
	&-\int_0^T e^{-r t}\left[\min _{\theta \in \Theta}\{(R'\nabla V(\tilde{Z}(t)) + v)\cdot\theta \} +\big(c-\Gamma\nabla V(\tilde{Z}(t)) \big)\cdot \tilde{Z}(t)\right] \mathrm{d} t\\
	=&\int_0^T e^{-r t} \nabla V(\tilde{Z}(t)) \cdot \mathrm{d} W(t) -\int_0^T e^{-r t} v_B \cdot \mathrm{d} \tilde{Y}_{B}(t)\\
    &-\int_0^T e^{-r t} F(\tilde{Z}(t), G(\tilde{Z}(t))) \mathrm{d} t.
	\end{aligned}
	$$

\end{proof}

\section{Model parameters of high-dimensional test problems}\label{app:hdtest}
   \paragraph{24-D (I) model}

The 24-D (I) model is constructed by concatenating three Zigzag models in the order Zigzag B, Zigzag C, and Zigzag A; see Figure~\ref{fig:24d}(a). We renumber the classes so that all left classes precede the right classes. To form a connected matching system, we add three basic activities linking classes 5–16, 9–20, and 1–24. To ensure these new activities are basic (i.e., receive nonzero resource allocation), we increase the arrival rates of the participating classes each by 1. We also renumber the activities by listing the basic activities first, followed by the nonbasic ones, so that the matrix $R$ and the matching rate $v$ can be properly partitioned. This setup ensures that Assumption~\ref{as:bf} is satisfied. In Figure~\ref{fig:24d}, only the numbering of basic activities is shown, corresponding to the column ordering of $R$.
   The arrival rates, matching values, and the optimal solution are given as follows:
$$
\begin{aligned}
&\lambda =  (4, 2, 2, 4, 4, 2, 2, 4, 4, 2, 1, 4, 4, 2, 2, 4, 4, 2, 2, 4, 4, 1, 2, 4)', \\
&v  = (2, 4, 8, 9, 8, 4, 2, 2, 4, 8, 9, 8, 4, 2, 2, 4, 8, 8, 4, 2, 4, 4, 4, 2, 4, 2, 6, 2, 2, 4, 2, 1)',\\
&x^* = (3, 1, 1, 1, 1, 1, 3, 3, 1, 1, 1, 1, 1, 3, 3, 1, 1, 1, 1, 3, 1, 1, 1, 0, 0, 0, 0, 0, 0, 0, 0, 0)'.
\end{aligned}
$$
For matching matrix $R_I$ see (\ref{eq:RI}).The abandonment and holding costs are uniformly sampled from $[0,10]$, and the abandonment rates are uniformly generated from $[0,0.4]$: 
$$\begin{aligned}
h  =& (5.6, 6.5, 4.8, 1.7, 0.0, 9.7, 5.7, 0.7, 4.0, 0.4, 7.7, 9.8, 1.0, 0.5, 0.3, 5.1, 7.8, 1.9, 1.9, 7.1, 9.8, 9.4, 6.4, 9.2)',\\
a = & (3.5, 9.1, 6.0, 5.5, 2.1, 2.4, 0.3, 7.8, 5.5, 3.7, 9.9, 4.6, 2.8, 8.1, 5.8, 8.2, 3.2, 4.5, 8.0, 4.6, 6.9, 7.6, 5.6, 6.5)',\\
\gamma = &(0.188, 0.291, 0.122, 0.355, 0.164, 0.287, 0.106, 0.098, 0.325, 0.199, 0.166, 0.291, 0.385, 0.124, 0.282, \\ &0.208, 0.293, 0.4,
0.083, 0.301, 0.187, 0.284, 0.349, 0.059)'.\\
\end{aligned}$$
  \paragraph{24-D (II) model}
  The 24-D (II) model is constructed by concatenating three copies of the most complex Zigzag model C; see Figure~\ref{fig:24d}(b). We add three activities linking classes 3–18, 7–22, and 11–14, thereby forming a connected network. The arrival rates of the participating classes are each increased by 1, but, according to the static planning problem, the resulting activities remain nonbasic. As with the 24-D (I) model, we renumber the classes so that left classes precede the right classes, and the activities are relabeled so that basic activities are listed before the nonbasic ones. The arrival rates, matching values, and the optimal solution are given as follows:
$$
\begin{aligned}
\lambda &= (3, 2, 3, 4, 3, 2, 3, 4, 3, 2, 3, 4, 4, 3, 2, 3, 4, 3, 2, 3, 4, 3, 2, 3)', \\
v  &= (2, 4, 8, 9, 8, 4, 2, 2, 4, 8, 9, 8, 4, 2, 2, 4, 8, 9, 8, 4, 2, 2, 2, 2, 2, 4, 2, 6, 2, 2, 4, 2, 2, 4, 2, 6, 2, 2, 4, \\&2, 2, 4, 2, 6, 2, 2, 4, 2)',\\
x^* &= (3, 1, 1, 2, 1, 1, 3, 3, 1, 1, 2, 1, 1, 3, 3, 1, 1, 2, 1, 1, 3, 0, 0, 0, 0, 0, 0, 0, 0, 0, 0, 0, 0, 0, 0, 0, 0, 0, 0, \\&0, 0, 0, 0, 0, 0, 0, 0, 0)'.
\end{aligned}
$$

For matching matrix $R_{II}$ see (\ref{eq:RII}). The abandonment and holding costs are uniformly sampled from $[0,20]$, and the abandonment rates are uniformly generated from $[0,0.2]$: 
$$\begin{aligned}
h  =& (11.2, 12.9, 9.6, 3.5, 0.0, 19.5, 11.4, 1.3, 8.0, 0.9, 15.3, 19.6, 2.0, 1.0, 0.6, 10.2, 15.6, 3.9, 3.8, 14.2, \\&19.6, 18.8, 12.8, 18.4)',
\end{aligned}$$
$$\begin{aligned}
a = & (7.1, 18.2, 11.9, 10.9, 4.2, 4.9, 0.6, 15.7, 11.1, 7.4, 19.8, 9.2, 5.6, 16.3, 11.5, 16.4, 6.5, 8.9, 16.0, 9.1, \\&13.8, 15.1, 11.3, 13.1)',\\
\gamma = &(0.094, 0.146, 0.061, 0.177, 0.082, 0.143, 0.053, 0.049, 0.163, 0.1, 0.083, 0.146, 0.193, 0.062, 0.141, 0.104,  \\
&0.146, 0.2, 0.041, 0.151, 0.094, 0.142, 0.174, 0.03)'.\\
\end{aligned}$$
\begin{equation}\label{eq:RI}
\setlength{\arraycolsep}{3pt}
\renewcommand{\arraystretch}{0.8} 
R_{\text{I}} = \left[\begin{array}{*{32}{c}}
1 & 0 & 0 & 0 & 0 & 0 & 0 & 0 & 0 & 0 & 0 & 0 & 0 & 0 & 0 & 0 & 0 & 0 & 0 & 0 & 0 & 0 & 1 & 0 & 0 & 0 & 0 & 0 & 0 & 0 & 0 & 0 \\
0 & 1 & 1 & 0 & 0 & 0 & 0 & 0 & 0 & 0 & 0 & 0 & 0 & 0 & 0 & 0 & 0 & 0 & 0 & 0 & 0 & 0 & 0 & 0 & 0 & 0 & 0 & 0 & 0 & 0 & 0 & 0 \\
0 & 0 & 0 & 1 & 1 & 0 & 0 & 0 & 0 & 0 & 0 & 0 & 0 & 0 & 0 & 0 & 0 & 0 & 0 & 0 & 0 & 0 & 0 & 0 & 0 & 0 & 0 & 0 & 0 & 0 & 0 & 0 \\
0 & 0 & 0 & 0 & 0 & 1 & 1 & 0 & 0 & 0 & 0 & 0 & 0 & 0 & 0 & 0 & 0 & 0 & 0 & 0 & 0 & 0 & 0 & 0 & 0 & 0 & 0 & 0 & 0 & 0 & 0 & 0 \\
0 & 0 & 0 & 0 & 0 & 0 & 0 & 1 & 0 & 0 & 0 & 0 & 0 & 0 & 0 & 0 & 0 & 0 & 0 & 0 & 1 & 0 & 0 & 1 & 1 & 1 & 0 & 0 & 0 & 0 & 0 & 0 \\
0 & 0 & 0 & 0 & 0 & 0 & 0 & 0 & 1 & 1 & 0 & 0 & 0 & 0 & 0 & 0 & 0 & 0 & 0 & 0 & 0 & 0 & 0 & 0 & 0 & 0 & 1 & 1 & 0 & 0 & 0 & 0 \\
0 & 0 & 0 & 0 & 0 & 0 & 0 & 0 & 0 & 0 & 1 & 1 & 0 & 0 & 0 & 0 & 0 & 0 & 0 & 0 & 0 & 0 & 0 & 0 & 0 & 0 & 0 & 0 & 1 & 1 & 0 & 0 \\
0 & 0 & 0 & 0 & 0 & 0 & 0 & 0 & 0 & 0 & 0 & 0 & 1 & 1 & 0 & 0 & 0 & 0 & 0 & 0 & 0 & 0 & 0 & 0 & 0 & 0 & 0 & 0 & 0 & 0 & 1 & 0 \\
0 & 0 & 0 & 0 & 0 & 0 & 0 & 0 & 0 & 0 & 0 & 0 & 0 & 0 & 1 & 0 & 0 & 0 & 0 & 0 & 0 & 1 & 0 & 0 & 0 & 0 & 0 & 0 & 0 & 0 & 0 & 0 \\
0 & 0 & 0 & 0 & 0 & 0 & 0 & 0 & 0 & 0 & 0 & 0 & 0 & 0 & 0 & 1 & 1 & 0 & 0 & 0 & 0 & 0 & 0 & 0 & 0 & 0 & 0 & 0 & 0 & 0 & 0 & 0 \\
0 & 0 & 0 & 0 & 0 & 0 & 0 & 0 & 0 & 0 & 0 & 0 & 0 & 0 & 0 & 0 & 0 & 1 & 0 & 0 & 0 & 0 & 0 & 0 & 0 & 0 & 0 & 0 & 0 & 0 & 0 & 1 \\
0 & 0 & 0 & 0 & 0 & 0 & 0 & 0 & 0 & 0 & 0 & 0 & 0 & 0 & 0 & 0 & 0 & 0 & 1 & 1 & 0 & 0 & 0 & 0 & 0 & 0 & 0 & 0 & 0 & 0 & 0 & 0 \\
1 & 1 & 0 & 0 & 0 & 0 & 0 & 0 & 0 & 0 & 0 & 0 & 0 & 0 & 0 & 0 & 0 & 0 & 0 & 0 & 0 & 0 & 0 & 0 & 0 & 0 & 0 & 0 & 0 & 0 & 0 & 0 \\
0 & 0 & 1 & 1 & 0 & 0 & 0 & 0 & 0 & 0 & 0 & 0 & 0 & 0 & 0 & 0 & 0 & 0 & 0 & 0 & 0 & 0 & 0 & 0 & 0 & 0 & 0 & 0 & 0 & 0 & 0 & 0 \\
0 & 0 & 0 & 0 & 1 & 1 & 0 & 0 & 0 & 0 & 0 & 0 & 0 & 0 & 0 & 0 & 0 & 0 & 0 & 0 & 0 & 0 & 0 & 0 & 0 & 0 & 0 & 0 & 0 & 0 & 0 & 0 \\
0 & 0 & 0 & 0 & 0 & 0 & 1 & 0 & 0 & 0 & 0 & 0 & 0 & 0 & 0 & 0 & 0 & 0 & 0 & 0 & 1 & 0 & 0 & 0 & 0 & 0 & 0 & 0 & 0 & 0 & 0 & 0 \\
0 & 0 & 0 & 0 & 0 & 0 & 0 & 1 & 1 & 0 & 0 & 0 & 0 & 0 & 0 & 0 & 0 & 0 & 0 & 0 & 0 & 0 & 0 & 0 & 0 & 0 & 0 & 0 & 1 & 0 & 0 & 0 \\
0 & 0 & 0 & 0 & 0 & 0 & 0 & 0 & 0 & 1 & 1 & 0 & 0 & 0 & 0 & 0 & 0 & 0 & 0 & 0 & 0 & 0 & 0 & 1 & 0 & 0 & 0 & 0 & 0 & 0 & 1 & 0 \\
0 & 0 & 0 & 0 & 0 & 0 & 0 & 0 & 0 & 0 & 0 & 1 & 1 & 0 & 0 & 0 & 0 & 0 & 0 & 0 & 0 & 0 & 0 & 0 & 1 & 0 & 1 & 0 & 0 & 0 & 0 & 0 \\
0 & 0 & 0 & 0 & 0 & 0 & 0 & 0 & 0 & 0 & 0 & 0 & 0 & 1 & 0 & 0 & 0 & 0 & 0 & 0 & 0 & 1 & 0 & 0 & 0 & 1 & 0 & 1 & 0 & 1 & 0 & 0 \\
0 & 0 & 0 & 0 & 0 & 0 & 0 & 0 & 0 & 0 & 0 & 0 & 0 & 0 & 1 & 1 & 0 & 0 & 0 & 0 & 0 & 0 & 0 & 0 & 0 & 0 & 0 & 0 & 0 & 0 & 0 & 0 \\
0 & 0 & 0 & 0 & 0 & 0 & 0 & 0 & 0 & 0 & 0 & 0 & 0 & 0 & 0 & 0 & 1 & 0 & 0 & 0 & 0 & 0 & 0 & 0 & 0 & 0 & 0 & 0 & 0 & 0 & 0 & 1 \\
0 & 0 & 0 & 0 & 0 & 0 & 0 & 0 & 0 & 0 & 0 & 0 & 0 & 0 & 0 & 0 & 0 & 1 & 1 & 0 & 0 & 0 & 0 & 0 & 0 & 0 & 0 & 0 & 0 & 0 & 0 & 0 \\
0 & 0 & 0 & 0 & 0 & 0 & 0 & 0 & 0 & 0 & 0 & 0 & 0 & 0 & 0 & 0 & 0 & 0 & 0 & 1 & 0 & 0 & 1 & 0 & 0 & 0 & 0 & 0 & 0 & 0 & 0 & 0
\end{array}\right]  
\end{equation}

\begin{equation}\label{eq:RII}
\setlength{\arraycolsep}{2pt}
R_{\text{II}} = 
\left[\begin{array}{*{48}{c}}
1 & 0 & 0 & 0 & 0 & 0 & 0 & 0 & 0 & 0 & 0 & 0 & 0 & 0 & 0 & 0 & 0 & 0 & 0 & 0 & 0 & 0 & 0 & 0 & 1 & 1 & 1 & 0 & 0 & 0 & 0 & 0 & 0 & 0 & 0 & 0 & 0 & 0 & 0 & 0 & 0 & 0 & 0 & 0 & 0 & 0 & 0 & 0 \\
0 & 1 & 1 & 0 & 0 & 0 & 0 & 0 & 0 & 0 & 0 & 0 & 0 & 0 & 0 & 0 & 0 & 0 & 0 & 0 & 0 & 0 & 0 & 0 & 0 & 0 & 0 & 1 & 1 & 0 & 0 & 0 & 0 & 0 & 0 & 0 & 0 & 0 & 0 & 0 & 0 & 0 & 0 & 0 & 0 & 0 & 0 & 0 \\
0 & 0 & 0 & 1 & 1 & 0 & 0 & 0 & 0 & 0 & 0 & 0 & 0 & 0 & 0 & 0 & 0 & 0 & 0 & 0 & 0 & 1 & 0 & 0 & 0 & 0 & 0 & 0 & 0 & 1 & 1 & 0 & 0 & 0 & 0 & 0 & 0 & 0 & 0 & 0 & 0 & 0 & 0 & 0 & 0 & 0 & 0 & 0 \\
0 & 0 & 0 & 0 & 0 & 1 & 1 & 0 & 0 & 0 & 0 & 0 & 0 & 0 & 0 & 0 & 0 & 0 & 0 & 0 & 0 & 0 & 0 & 0 & 0 & 0 & 0 & 0 & 0 & 0 & 0 & 1 & 0 & 0 & 0 & 0 & 0 & 0 & 0 & 0 & 0 & 0 & 0 & 0 & 0 & 0 & 0 & 0\\
0 & 0 & 0 & 0 & 0 & 0 & 0 & 1 & 0 & 0 & 0 & 0 & 0 & 0 & 0 & 0 & 0 & 0 & 0 & 0 & 0 & 0 & 0 & 0 & 0 & 0 & 0 & 0 & 0 & 0 & 0 & 0 & 1 & 1 & 1 & 0 & 0 & 0 & 0 & 0 & 0 & 0 & 0 & 0 & 0 & 0 & 0 & 0 \\
0 & 0 & 0 & 0 & 0 & 0 & 0 & 0 & 1 & 1 & 0 & 0 & 0 & 0 & 0 & 0 & 0 & 0 & 0 & 0 & 0 & 0 & 0 & 0 & 0 & 0 & 0 & 0 & 0 & 0 & 0 & 0 & 0 & 0 & 0 & 1 & 1 & 0 & 0 & 0 & 0 & 0 & 0 & 0 & 0 & 0 & 0 & 0 \\
0 & 0 & 0 & 0 & 0 & 0 & 0 & 0 & 0 & 0 & 1 & 1 & 0 & 0 & 0 & 0 & 0 & 0 & 0 & 0 & 0 & 0 & 1 & 0 & 0 & 0 & 0 & 0 & 0 & 0 & 0 & 0 & 0 & 0 & 0 & 0 & 0 & 1 & 1 & 0 & 0 & 0 & 0 & 0 & 0 & 0 & 0 & 0 \\
0 & 0 & 0 & 0 & 0 & 0 & 0 & 0 & 0 & 0 & 0 & 0 & 1 & 1 & 0 & 0 & 0 & 0 & 0 & 0 & 0 & 0 & 0 & 0 & 0 & 0 & 0 & 0 & 0 & 0 & 0 & 0 & 0 & 0 & 0 & 0 & 0 & 0 & 0 & 1 & 0 & 0 & 0 & 0 & 0 & 0 & 0 & 0\\
0 & 0 & 0 & 0 & 0 & 0 & 0 & 0 & 0 & 0 & 0 & 0 & 0 & 0 & 1 & 0 & 0 & 0 & 0 & 0 & 0 & 0 & 0 & 0 & 0 & 0 & 0 & 0 & 0 & 0 & 0 & 0 & 0 & 0 & 0 & 0 & 0 & 0 & 0 & 0 & 1 & 1 & 1 & 0 & 0 & 0 & 0 & 0 \\
0 & 0 & 0 & 0 & 0 & 0 & 0 & 0 & 0 & 0 & 0 & 0 & 0 & 0 & 0 & 1 & 1 & 0 & 0 & 0 & 0 & 0 & 0 & 0 & 0 & 0 & 0 & 0 & 0 & 0 & 0 & 0 & 0 & 0 & 0 & 0 & 0 & 0 & 0 & 0 & 0 & 0 & 0 & 1 & 1 & 0 & 0 & 0 \\
0 & 0 & 0 & 0 & 0 & 0 & 0 & 0 & 0 & 0 & 0 & 0 & 0 & 0 & 0 & 0 & 0 & 1 & 1 & 0 & 0 & 0 & 0 & 1 & 0 & 0 & 0 & 0 & 0 & 0 & 0 & 0 & 0 & 0 & 0 & 0 & 0 & 0 & 0 & 0 & 0 & 0 & 0 & 0 & 0 & 1 & 1 & 0 \\
0 & 0 & 0 & 0 & 0 & 0 & 0 & 0 & 0 & 0 & 0 & 0 & 0 & 0 & 0 & 0 & 0 & 0 & 0 & 1 & 1 & 0 & 0 & 0 & 0 & 0 & 0 & 0 & 0 & 0 & 0 & 0 & 0 & 0 & 0 & 0 & 0 & 0 & 0 & 0 & 0 & 0 & 0 & 0 & 0 & 0 & 0 & 1\\
1 & 1 & 0 & 0 & 0 & 0 & 0 & 0 & 0 & 0 & 0 & 0 & 0 & 0 & 0 & 0 & 0 & 0 & 0 & 0 & 0 & 0 & 0 & 0 & 0 & 0 & 0 & 0 & 0 & 1 & 0 & 0 & 0 & 0 & 0 & 0 & 0 & 0 & 0 & 0 & 0 & 0 & 0 & 0 & 0 & 0 & 0 & 0 \\
0 & 0 & 1 & 1 & 0 & 0 & 0 & 0 & 0 & 0 & 0 & 0 & 0 & 0 & 0 & 0 & 0 & 0 & 0 & 0 & 0 & 0 & 0 & 1 & 1 & 0 & 0 & 0 & 0 & 0 & 0 & 1 & 0 & 0 & 0 & 0 & 0 & 0 & 0 & 0 & 0 & 0 & 0 & 0 & 0 & 0 & 0 & 0 \\
0 & 0 & 0 & 0 & 1 & 1 & 0 & 0 & 0 & 0 & 0 & 0 & 0 & 0 & 0 & 0 & 0 & 0 & 0 & 0 & 0 & 0 & 0 & 0 & 0 & 1 & 0 & 1 & 0 & 0 & 0 & 0 & 0 & 0 & 0 & 0 & 0 & 0 & 0 & 0 & 0 & 0 & 0 & 0 & 0 & 0 & 0 & 0 \\
0 & 0 & 0 & 0 & 0 & 0 & 1 & 0 & 0 & 0 & 0 & 0 & 0 & 0 & 0 & 0 & 0 & 0 & 0 & 0 & 0 & 0 & 0 & 0 & 0 & 0 & 1 & 0 & 1 & 0 & 1 & 0 & 0 & 0 & 0 & 0 & 0 & 0 & 0 & 0 & 0 & 0 & 0 & 0 & 0 & 0 & 0 & 0\\
0 & 0 & 0 & 0 & 0 & 0 & 0 & 1 & 1 & 0 & 0 & 0 & 0 & 0 & 0 & 0 & 0 & 0 & 0 & 0 & 0 & 0 & 0 & 0 & 0 & 0 & 0 & 0 & 0 & 0 & 0 & 0 & 0 & 0 & 0 & 0 & 0 & 1 & 0 & 0 & 0 & 0 & 0 & 0 & 0 & 0 & 0 & 0 \\
0 & 0 & 0 & 0 & 0 & 0 & 0 & 0 & 0 & 1 & 1 & 0 & 0 & 0 & 0 & 0 & 0 & 0 & 0 & 0 & 0 & 1 & 0 & 0 & 0 & 0 & 0 & 0 & 0 & 0 & 0 & 0 & 1 & 0 & 0 & 0 & 0 & 0 & 0 & 1 & 0 & 0 & 0 & 0 & 0 & 0 & 0 & 0 \\
0 & 0 & 0 & 0 & 0 & 0 & 0 & 0 & 0 & 0 & 0 & 1 & 1 & 0 & 0 & 0 & 0 & 0 & 0 & 0 & 0 & 0 & 0 & 0 & 0 & 0 & 0 & 0 & 0 & 0 & 0 & 0 & 0 & 1 & 0 & 1 & 0 & 0 & 0 & 0 & 0 & 0 & 0 & 0 & 0 & 0 & 0 & 0 \\
0 & 0 & 0 & 0 & 0 & 0 & 0 & 0 & 0 & 0 & 0 & 0 & 0 & 1 & 0 & 0 & 0 & 0 & 0 & 0 & 0 & 0 & 0 & 0 & 0 & 0 & 0 & 0 & 0 & 0 & 0 & 0 & 0 & 0 & 1 & 0 & 1 & 0 & 1 & 0 & 0 & 0 & 0 & 0 & 0 & 0 & 0 & 0\\
0 & 0 & 0 & 0 & 0 & 0 & 0 & 0 & 0 & 0 & 0 & 0 & 0 & 0 & 1 & 1 & 0 & 0 & 0 & 0 & 0 & 0 & 0 & 0 & 0 & 0 & 0 & 0 & 0 & 0 & 0 & 0 & 0 & 0 & 0 & 0 & 0 & 0 & 0 & 0 & 0 & 0 & 0 & 0 & 0 & 1 & 0 & 0 \\
0 & 0 & 0 & 0 & 0 & 0 & 0 & 0 & 0 & 0 & 0 & 0 & 0 & 0 & 0 & 0 & 1 & 1 & 0 & 0 & 0 & 0 & 1 & 0 & 0 & 0 & 0 & 0 & 0 & 0 & 0 & 0 & 0 & 0 & 0 & 0 & 0 & 0 & 0 & 0 & 1 & 0 & 0 & 0 & 0 & 0 & 0 & 1 \\
0 & 0 & 0 & 0 & 0 & 0 & 0 & 0 & 0 & 0 & 0 & 0 & 0 & 0 & 0 & 0 & 0 & 0 & 1 & 1 & 0 & 0 & 0 & 0 & 0 & 0 & 0 & 0 & 0 & 0 & 0 & 0 & 0 & 0 & 0 & 0 & 0 & 0 & 0 & 0 & 0 & 1 & 0 & 1 & 0 & 0 & 0 & 0 \\
0 & 0 & 0 & 0 & 0 & 0 & 0 & 0 & 0 & 0 & 0 & 0 & 0 & 0 & 0 & 0 & 0 & 0 & 0 & 0 & 1 & 0 & 0 & 0 & 0 & 0 & 0 & 0 & 0 & 0 & 0 & 0 & 0 & 0 & 0 & 0 & 0 & 0 & 0 & 0 & 0 & 0 & 1 & 0 & 1 & 0 & 1 & 0
\end{array}\right]
\end{equation}

\begin{figure}[H]
    \centering
    \begin{subfigure}[b]{0.4\textwidth}
        \centering
    \includegraphics[width=0.7\textwidth]{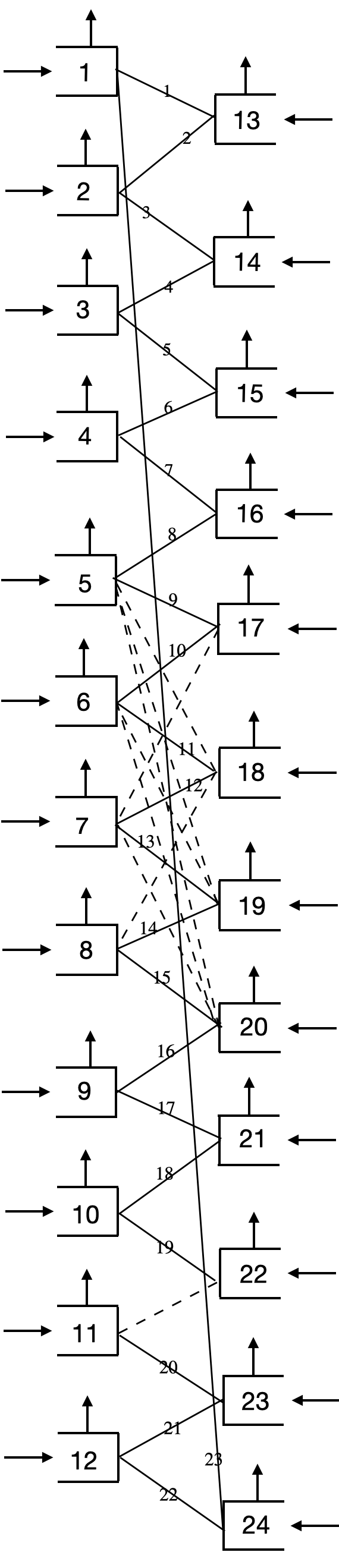}
        \subcaption{24-D (I)}
        \label{fig:24d_1}
    \end{subfigure}
    \hspace{0.04\textwidth}
    \begin{subfigure}[b]{0.4\textwidth}
        \centering
        \includegraphics[width=0.7\textwidth]{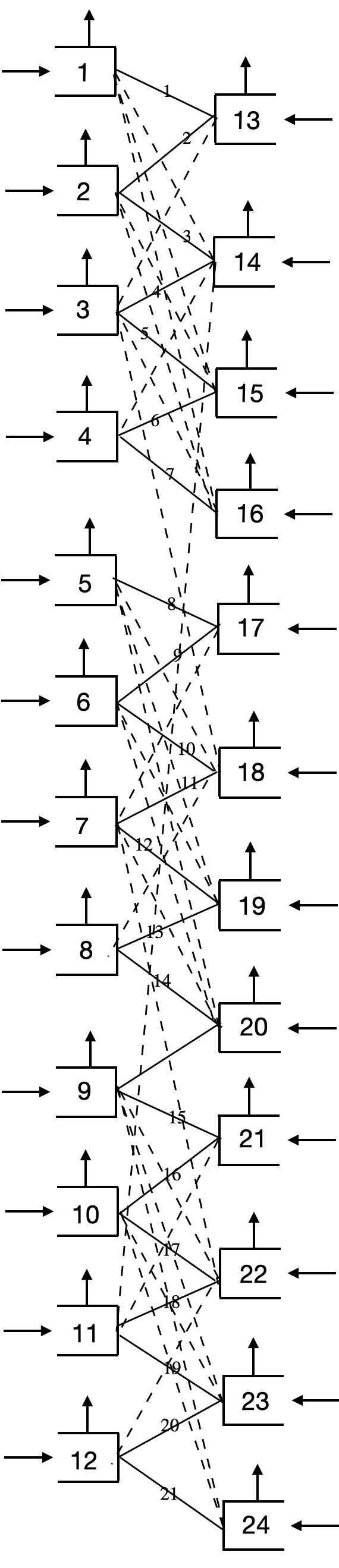}
        \subcaption{24-D (II)}
        \label{fig:24d_2}
    \end{subfigure}
    \caption{24-dimensional models.}
    \label{fig:24d}
\end{figure}

   \paragraph{120-D model} The 120-D model is constructed by concatenating five copies of the 24-D (II) model. We renumber the classes so that left classes precede the right classes. To ensure connectivity, we add four basic activities linking classes 12–73, 24–85, 36–97, and 48–109. To make these added activities basic, we increase the arrival rates of the eight participating classes by 1, as in the 24-D (I) model. The activities are then renumbered so that basic activities are listed before the nonbasic ones. Due to space constraints, the figure of the 120-D model is omitted, but its structure can be fully reconstructed from the description and the corresponding $R$ matrix. The arrival rates, matching values, and corresponding optimal solution are given as follows:
$$
\begin{aligned}
\lambda = & (3, 2, 3, 4, 3, 2, 3, 4, 3, 2, 3, 5, 3, 2, 3, 4, 3, 2, 3, 4, 3, 2, 3, 5, 3, 2, 3, 4, 3, 2, 3, 4, 3, 2, 3, 5, 3, 2, 3, 4, 3, 2,\\& 3, 4, 3, 2, 3, 5, 3, 2, 3, 4, 3, 2, 3, 4, 3, 2, 3, 4, 4, 3, 2, 3, 4, 3, 2, 3, 4, 3, 2, 3, 5, 3, 2, 3, 4, 3, 2, 3, 4, 3, 2, 3, \\&5, 3, 2, 3, 4, 3, 2, 3, 4, 3, 2, 3, 5, 3, 2, 3, 4, 3, 2, 3, 4, 3, 2, 3, 5, 3, 2, 3, 4, 3, 2, 3, 4, 3, 2, 3)',\\
v = & (2, 4, 8, 9, 8, 4, 2, 2, 4, 8, 9, 8, 4, 2, 2, 4, 8, 9, 8, 4, 2, 2, 4, 8, 9, 8, 4, 2, 2, 4, 8, 9, 8, 4, 2, 2, 4, 8, 9, 8, 4, 2, \\&2, 4, 8, 9, 8, 4, 2, 2, 4, 8, 9, 8, 4, 2, 2, 4, 8, 9, 8, 4, 2, 2, 4, 8, 9, 8, 4, 2, 2, 4, 8, 9, 8, 4, 2, 2, 4, 8, 9, 8, 4, 2, \\& 2, 4, 8, 9, 8, 4, 2, 2, 4, 8, 9, 8, 4, 2, 2, 4, 8, 9, 8, 4, 2, 1, 1, 1, 1, 2, 2, 2, 2, 4, 2, 6, 2, 2, 4, 2, 2, 4, 2, 6, 2, 2, \\& 4,  2, 2, 4, 2, 6, 2,  2, 4, 2, 2, 2, 2, 2, 4, 2, 6, 2, 2, 4, 2, 2, 4, 2, 6, 2, 2, 4, 2, 2, 4, 2, 6, 2, 2, 4, 2, 2, 2, 2, 2, 4,\\& 2, 6, 2, 2, 4, 2, 2, 4, 2, 6, 2, 2, 4, 2, 2, 4, 2, 6, 2, 2, 4, 2, 2, 2, 2, 2, 4, 2, 6, 2, 2, 4, 2, 2, 4, 2, 6, 2, 2, 4, 2, 2, \\& 4, 2, 6, 2, 2, 4, 2, 2, 2, 2, 2, 4, 2, 6, 2, 2, 4, 2, 2, 4, 2, 6, 2, 2, 4, 2, 2, 4, 2, 6, 2, 2, 4, 2)',\\
x^* = & (3, 1, 1, 2, 1, 1, 3, 3, 1, 1, 2, 1, 1, 3, 3, 1, 1, 2, 1, 1, 3, 3, 1, 1, 2, 1, 1, 3, 3, 1, 1, 2, 1, 1, 3, 3, 1, 1, 2, 1, 1, 3, \\&3,  1, 1, 2, 1, 1, 3, 3, 1, 1, 2, 1, 1, 3, 3, 1, 1, 2, 1, 1, 3, 3, 1, 1, 2, 1, 1, 3, 3, 1, 1, 2, 1, 1, 3, 3, 1, 1, 2, 1, 1, 3, 3,\\& 1, 1, 2, 1, 1, 3, 3, 1, 1, 2, 1, 1, 3, 3, 1, 1, 2, 1, 1, 3, 1, 1, 1, 1, 0, 0, 0, 0, 0, 0, 0, 0, 0, 0, 0, 0, 0, 0, 0, 0, 0, 0, 0,\\& 0,  0, 0, 0, 0, 0, 0, 0, 0, 0, 0, 0, 0, 0, 0, 0, 0, 0, 0, 0, 0, 0, 0, 0, 0, 0, 0, 0, 0, 0, 0, 0, 0, 0, 0, 0, 0, 0, 0, 0, 0, 0, 0,\\& 0, 0, 0, 0, 0, 0, 0, 0, 0, 0, 0, 0, 0, 0, 0, 0, 0, 0, 0, 0, 0, 0, 0, 0, 0, 0, 0, 0, 0, 0, 0, 0, 0, 0, 0, 0, 0, 0, 0, 0, 0, 0, 0,\\& 0, 0, 0, 0, 0, 0, 0, 0, 0, 0, 0, 0, 0, 0, 0, 0, 0, 0, 0, 0, 0, 0, 0, 0, 0, 0, 0, 0, 0, 0)'.
\end{aligned}
$$The abandonment and holding costs are uniformly sampled from $[0,10]$, and the abandonment rates are uniformly generated from $[0,1]$:
$$\begin{aligned}
h = &(5.6, 6.5, 4.8, 1.7, 0.0, 9.7, 5.7, 0.7, 4.0, 0.4, 7.7, 9.8, 1.0, 0.5, 0.3, 5.1, 7.8, 1.9, 1.9, 7.1, 9.8, 9.4, 6.4, 9.2,\\& 3.5, 9.1, 6.0, 5.5, 2.1, 2.4, 0.3, 7.8, 5.5, 3.7, 9.9, 4.6, 2.8, 8.1, 5.8, 8.2, 3.2, 4.5, 8.0, 4.6, 6.9, 7.6, 5.6, 6.5, \\& 5.9, 7.6, 8.6, 8.4, 3.4, 5.4,  6.0, 1.6, 2.8, 3.3, 3.6, 5.6, 6.7, 2.8, 2.7, 6.4, 1.5, 3.5, 7.3, 6.2, 8.1, 5.4, 5.7, 5.0, \\&3.8, 3.4, 2.5, 5.9, 5.1, 5.7, 6.9, 2.1, 2.0, 1.4, 4.4, 5.7, 7.6, 2.6, 0.1, 5.8, 2.0, 1.4, 1.7, 6.8, 7.3, 2.6, 9.1, 1.4, \\& 6.3, 7.4, 1.5, 9.5, 2.3, 7.5, 7.9, 0.8, 5.1, 2.4, 4.2, 8.0,  8.0, 6.8, 5.9, 2.3, 8.9, 3.5, 7.1, 5.2, 3.7, 4.8, 2.6, 0.7)',\\\end{aligned}$$
$$\begin{aligned}a = &(9.3, 3.1, 5.1, 6.2, 6.2, 4.2, 3.3, 7.8, 9.6, 9.6, 9.1, 1.6, 8.6, 3.3, 7.4, 5.0, 1.8, 5.1, 0.7, 1.8, 3.5, 7.9, 0.6, 8.3, \\&9.3, 1.2, 6.9, 6.0, 3.4, 3.0, 6.8, 6.1, 1.9, 0.5, 4.0, 5.8, 7.2, 2.8, 2.7, 3.2, 9.0, 4.6, 9.0, 10.0, 0.5, 0.2, 5.3, 8.2, \\& 7.1, 4.7, 1.5, 8.8, 6.3, 8.3, 8.3, 4.0, 3.5, 0.5, 0.9, 8.4, 2.9, 8.1, 9.1, 9.3, 4.3, 5.4, 9.0, 2.5, 3.4, 3.7, 7.6, 0.3, \\&8.4, 5.4, 4.6, 3.1, 3.4, 2.6, 6.7, 0.0, 6.2, 3.2, 6.5, 0.5, 6.7, 9.7, 3.0, 6.6, 0.4, 8.5, 3.7, 9.9, 3.3, 3.5, 1.3, 1.4, \\&6.1,  2.3, 0.3, 5.0, 4.8, 9.7, 2.0, 6.1, 9.9, 4.2, 9.1, 1.5, 6.9, 8.2, 2.4, 4.2, 5.3, 1.3, 6.5, 7.0, 2.1, 5.2, 8.2, 3.2)',\\\end{aligned}$$
$$\begin{aligned}\gamma = &(0.47, 0.728, 0.304, 0.887, 0.41, 0.717, 0.265, 0.245, 0.813, 0.498, 0.416, 0.728, 0.963, 0.31, 0.704, \\&0.519, 0.731, 1.0,  0.206, 0.753, 0.469, 0.709, 0.872, 0.148, 0.213, 0.412, 0.058, 0.349, 0.417, 0.124,\\& 0.744, 0.763, 0.39, 0.345, 0.201,  0.427, 0.316, 0.214, 0.868, 0.229, 0.04, 0.225, 0.019, 0.865, 0.844,\\& 0.319, 0.96, 0.804, 0.421, 0.112, 0.851, 0.607, 0.231, 0.995, 0.366, 0.203, 0.493, 0.837, 0.141, 0.387,\\& 0.337, 0.94, 0.999, 0.465, 0.17, 0.696, 0.861, 0.332, 0.207, 0.687, 0.098, 0.844, 0.004, 0.151, 0.704, \\&0.35, 0.081, 0.803, 0.239, 0.461, 0.264, 0.524, 0.436, 0.972, 0.201, 0.071, 0.303, 0.136, 0.662,  0.25, \\& 0.101, 0.219, 0.195, 0.389, 0.515, 0.223, 0.355, 0.703, 0.8, 0.579, 0.932,0.544, 0.936, 0.72, 0.633, \\&0.119, 0.056, 0.125, 0.83, 0.907,  0.615, 0.088, 0.504, 0.236, 0.584, 0.732, 0.162, 0.125, 0.104, 0.944)'.
\end{aligned}$$
For the 120-dimensional model, the matrix $R$ is too large to display in full. Each column of $R$ contains exactly two entries equal to 1 (all others are 0), indicating the classes participating in the corresponding activity. Instead of showing the complete matrix, we list the pairs of class indices corresponding to the two 1s in each column:
\{(1, 61), (2, 61), (2, 62), (3, 62), (3, 63), (4, 63), (4, 64), (5, 65), (6, 65), (6, 66), (7, 66), (7, 67), (8, 67), (8, 68), (9, 69), (10, 69), (10, 70), (11, 70), (11, 71), (12, 71), (12, 72), (13, 73), (14, 73), (14, 74), (15, 74), (15, 75), (16, 75), (16, 76), (17, 77), (18, 77), (18, 78), (19, 78), (19, 79), (20, 79), (20, 80), (21, 81), (22, 81), (22, 82), (23, 82), (23, 83), (24, 83), (24, 84), (25, 85), (26, 85), (26, 86), (27, 86), (27, 87), (28, 87), (28, 88), (29, 89), (30, 89), (30, 90), (31, 90), (31, 91), (32, 91), (32, 92), (33, 93), (34, 93), (34, 94), (35, 94), (35, 95), (36, 95), (36, 96), (37, 97), (38, 97), (38, 98), (39, 98), (39, 99), (40, 99), (40, 100), (41, 101), (42, 101), (42, 102), (43, 102), (43, 103), (44, 103), (44, 104), (45, 105), (46, 105), (46, 106), (47, 106), (47, 107), (48, 107), (48, 108), (49, 109), (50, 109), (50, 110), (51, 110), (51, 111), (52, 111), (52, 112), (53, 113), (54, 113), (54, 114), (55, 114), (55, 115), (56, 115), (56, 116), (57, 117), (58, 117), (58, 118), (59, 118), (59, 119), (60, 119), (60, 120), (12, 73), (24, 85), (36, 97), (48, 109), (3, 66), (7, 70), (11, 62), (1, 62), (1, 63), (1, 64), (2, 63), (2, 64), (3, 61), (3, 64), (4, 62), (5, 66), (5, 67), (5, 68), (6, 67), (6, 68), (7, 65), (7, 68), (8, 66), (9, 70), (9, 71), (9, 72), (10, 71), (10, 72), (11, 69), (11, 72), (12, 70), (15, 78), (19, 82), (23, 74), (13, 74), (13, 75), (13, 76), (14, 75), (14, 76), (15, 73), (15, 76), (16, 74), (17, 78), (17, 79), (17, 80), (18, 79), (18, 80), (19, 77), (19, 80), (20, 78), (21, 82), (21, 83), (21, 84), (22, 83), (22, 84), (23, 81), (23, 84), (24, 82), (27, 90), (31, 94), (35, 86), (25, 86), (25, 87), (25, 88), (26, 87), (26, 88), (27, 85), (27, 88), (28, 86), (29, 90), (29, 91), (29, 92), (30, 91), (30, 92), (31, 89), (31, 92), (32, 90), (33, 94), (33, 95), (33, 96), (34, 95), (34, 96), (35, 93), (35, 96), (36, 94), (39, 102), (43, 106), (47, 98), (37, 98), (37, 99), (37, 100), (38, 99), (38, 100), (39, 97), (39, 100), (40, 98), (41, 102), (41, 103), (41, 104), (42, 103), (42, 104), (43, 101), (43, 104), (44, 102), (45, 106), (45, 107), (45, 108), (46, 107), (46, 108), (47, 105), (47, 108), (48, 106), (51, 114), (55, 118), (59, 110), (49, 110), (49, 111), (49, 112), (50, 111), (50, 112), (51, 109), (51, 112), (52, 110), (53, 114), (53, 115), (53, 116), (54, 115), (54, 116), (55, 113), (55, 116), (56, 114), (57, 118), (57, 119), (57, 120), (58, 119), (58, 120), (59, 117), (59, 120), (60, 118)\}.

   \section{Implementation Details of Our Computational Method}\label{app:NNarchi}
This section provides an overview of the computing environment and the neural network training setup for our nine test problems. We implement our method using fully connected deep neural networks in TensorFlow. Each model consists of three hidden layers with 100 neurons per layer and uses the elu activation function. Training is carried out with the Adam optimizer, a batch size of 256, and a stepwise learning rate schedule decaying from $10^{-3}$ to $5\times 10^{-6}$ over 80,000 iterations. Details, including problem-specific parameters, are summarized in Tables~\ref{table:parX}-\ref{table:parHD}.

All neural networks were trained on two high-performance clusters. The first cluster provided GPU nodes with NVIDIA V100 (16 GB), RTX 6000 (24 GB), and A100 (40 GB) accelerators, paired with Intel Xeon Gold 6248R CPUs. The second cluster consisted of nodes with 8× NVIDIA H100 (80 GB) GPUs per node, supported by 64-core CPUs and approximately 2 TB of system memory. The latter setup improved performance, reducing runtime by about 40\% relative to the former.

\begin{table}[H]
\renewcommand{\arraystretch}{1.0}
\centering
\begin{tabular}{@{}lccc@{}}
\toprule
\textbf{Hyperparameters} & \textbf{X High} & \textbf{X Medium} & \textbf{X Low} \\
\midrule
Horizon $T$ & 0.1 & 0.1 & 0.1 \\
Number of discretization intervals $N$ & 64 & 64 & 64 \\
Number of hidden layers & 3 & 3 & 3 \\
Number of neurons per layer & 100 & 100 & 100 \\
Activation function & \texttt{elu} & \texttt{elu} & \texttt{elu} \\
Optimizer & \texttt{Adam} & \texttt{Adam} & \texttt{Adam} \\
Batch size & 256 & 256 & 256 \\
\midrule
Drift bound $\eta$ & 0.5 & 0.5 & 1 \\
Rounding $\varepsilon$ & 0.03 & 0.03 & 0.15 \\
Reference policy $\tilde{\theta}$ (basic) & 0.1 & 0.1 & 0.1 \\
Reference policy $\tilde{\theta}$ (nonbasic) & $-0.01$ & $-0.01$ & $-0.01$ \\
\midrule
Number of iterations & 80000 & 80000 & 80000 \\
Learning rate (iteration range) &
\makecell[l]{$10^{-3}$ (0, $2\mathrm{e}4$), \\
$10^{-4}$ ($2\mathrm{e}4$, $5\mathrm{e}4$),\\ 
$10^{-5}$ ($5\mathrm{e}4$, $7\mathrm{e}4$), \\
$5 \times 10^{-6}$ ($7\mathrm{e}4$, $8\mathrm{e}4$)} &
\makecell[l]{$10^{-3}$ (0, $2\mathrm{e}4$), \\
$10^{-4}$ ($2\mathrm{e}4$, $5\mathrm{e}4$),\\ 
$10^{-5}$ ($5\mathrm{e}4$, $7\mathrm{e}4$), \\
$5 \times 10^{-6}$ ($7\mathrm{e}4$, $8\mathrm{e}4$)} &
\makecell[l]{$10^{-3}$ (0, $2\mathrm{e}4$), \\
$10^{-4}$ ($2\mathrm{e}4$, $5\mathrm{e}4$),\\ 
$10^{-5}$ ($5\mathrm{e}4$, $7\mathrm{e}4$), \\
$5 \times 10^{-6}$ ($7\mathrm{e}4$, $8\mathrm{e}4$)} \\
\bottomrule
\end{tabular}
\caption{Summary of the hyperparameters used for training X models.}
\label{table:parX}
\end{table}

\begin{table}[H]
\renewcommand{\arraystretch}{1.0}
\centering
\begin{tabular}{@{}lccc@{}}
\toprule
\textbf{Hyperparameters} & \textbf{Zigzag A} & \textbf{Zigzag B} & \textbf{Zigzag C} \\
\midrule
Horizon $T$ & 0.1 & 0.1 & 0.1 \\
Number of discretization intervals $N$ & 64 & 64 & 64 \\
Number of hidden layers & 3 & 3 & 3 \\
Number of neurons per layer & 100 & 100 & 100 \\
Activation function & \texttt{elu} & \texttt{elu} & \texttt{elu} \\
Optimizer & \texttt{Adam} & \texttt{Adam} & \texttt{Adam} \\
Batch size & 256 & 256 & 256 \\
\midrule
Drift bound $\eta$ & 1 & 0.5 & 1 \\
Rounding $\varepsilon$ & 0.05 & 0.07 & 0.05 \\
Reference policy $\tilde{\theta}$ (basic) & 3.0 & 3.0 & 2.0 \\
Reference policy $\tilde{\theta}$ (nonbasic) & $-0.1$ & -- & $-0.01$ \\
\midrule
Number of iterations & 80000 & 80000 & 80000 \\
Learning rate (iteration range) &
\makecell[l]{$10^{-3}$ (0, $2\mathrm{e}4$), \\
$10^{-4}$ ($2\mathrm{e}4$, $5\mathrm{e}4$),\\ 
$10^{-5}$ ($5\mathrm{e}4$, $7\mathrm{e}4$), \\
$5 \times 10^{-6}$ ($7\mathrm{e}4$, $8\mathrm{e}4$)} &
\makecell[l]{$10^{-3}$ (0, $2\mathrm{e}4$), \\
$10^{-4}$ ($2\mathrm{e}4$, $5\mathrm{e}4$),\\ 
$10^{-5}$ ($5\mathrm{e}4$, $7\mathrm{e}4$), \\
$5 \times 10^{-6}$ ($7\mathrm{e}4$, $8\mathrm{e}4$)} &
\makecell[l]{$10^{-3}$ (0, $2\mathrm{e}4$), \\
$10^{-4}$ ($2\mathrm{e}4$, $5\mathrm{e}4$),\\ 
$10^{-5}$ ($5\mathrm{e}4$, $7\mathrm{e}4$), \\
$5 \times 10^{-6}$ ($7\mathrm{e}4$, $8\mathrm{e}4$)} \\
\bottomrule
\end{tabular}
\caption{Summary of the hyperparameters used for training Zigzag models.}
\label{table:parZZ}
\end{table}

\begin{table}[H]
\renewcommand{\arraystretch}{1.0}
\centering
\begin{tabular}{@{}lccc@{}}
\toprule
\textbf{Hyperparameters} & \textbf{24-D (I)} & \textbf{24-D (II)} & \textbf{120-D} \\
\midrule
Horizon $T$ & 0.1 & 0.1 & 0.1 \\
Number of discretization intervals $N$ & 64 & 64 & 64 \\
Number of hidden layers & 3 & 3 & 3 \\
Number of neurons per layer & 100 & 100 & 100 \\
Activation function & \texttt{elu} & \texttt{elu} & \texttt{elu} \\
Optimizer & \texttt{Adam} & \texttt{Adam} & \texttt{Adam} \\
Batch size & 256 & 256 & 256 \\
\midrule
Drift bound $\eta$ & 0.5 & 1 & 0.0 \\
Rounding $\varepsilon$ & 0.03 & 0.0 & 0.0 \\
Reference policy $\tilde{\theta}$ (basic) & 0.1 & 0.1 & 0.1 \\
Reference policy $\tilde{\theta}$ (nonbasic) & $-0.1$ & $-0.01$ & $-0.01$ \\
\midrule
Number of iterations & 80000 & 80000 & 80000 \\
Learning rate (iteration range) &
\makecell[l]{$10^{-3}$ (0, $2\mathrm{e}4$), \\
$10^{-4}$ ($2\mathrm{e}4$, $5\mathrm{e}4$),\\ 
$10^{-5}$ ($5\mathrm{e}4$, $7\mathrm{e}4$), \\
$5 \times 10^{-6}$ ($7\mathrm{e}4$, $8\mathrm{e}4$)} &
\makecell[l]{$10^{-3}$ (0, $2\mathrm{e}4$), \\
$10^{-4}$ ($2\mathrm{e}4$, $5\mathrm{e}4$),\\ 
$10^{-5}$ ($5\mathrm{e}4$, $7\mathrm{e}4$), \\
$5 \times 10^{-6}$ ($7\mathrm{e}4$, $8\mathrm{e}4$)} &
\makecell[l]{$10^{-3}$ (0, $2\mathrm{e}4$), \\
$10^{-4}$ ($2\mathrm{e}4$, $5\mathrm{e}4$),\\ 
$10^{-5}$ ($5\mathrm{e}4$, $7\mathrm{e}4$), \\
$5 \times 10^{-6}$ ($7\mathrm{e}4$, $8\mathrm{e}4$)} \\
\bottomrule
\end{tabular}
\caption{Summary of the hyperparameters used for training high-dimensional models.}
\label{table:parHD}
\end{table}
Since the static-priority policy is of the discrete-review type, we determine the optimal review period $l$ for each test problem by performing a bisection search, followed by rounding to the nearest multiple of 0.005 or 0.01 to ensure consistency across experiments. Table~\ref{table:review_period} summarizes those best performing review periods for our nine test problems. 
\begin{table}[h]
		\centering
		\begin{tabular}{ccccccc}
	\hline \text { Model } & \text { Best review period $l$} \\
	\hline 
	\text { X High } & $0.001$\\
    \text { X Medium } & $0.001$\\
    \text { X Low } & $0.001$\\
	\text { Zigzag A} & $0.0001$ \\
	\text { Zigzag B} & $0.0001$ \\
	\text { Zigzag C} & $0.01$\\
	\text { 24-D (I)} & $0.0001$\\
    \text { 24-D (II)} & $0.005$\\
	\text { 120-D } & $0.01$\\
	\hline
		\end{tabular}
		\caption{Best performing discrete review periods $l$ for the static-priority policy.}
		\label{table:review_period}
	\end{table}

	\end{appendix}

\end{document}